\documentclass[12pt]{elsarticle}
\usepackage{mystyleElsArticle}
\usepackage{mathrsfs}
\usepackage{algorithm} 
\usepackage{booktabs}
\usepackage{bm}

\newcommand{\gvec}[1]{\boldsymbol{#1}}
\newcommand{\vxi}{\gvec{\xi}}

\newcommand{\mset}[1]{\mathfrak{#1}}

\usepackage[margin=1.00 in]{geometry}    

\makeatletter
\def\ps@pprintTitle{%
  \let\@oddhead\@empty
  \let\@evenhead\@empty
  \let\@oddfoot\@empty
  \let\@evenfoot\@oddfoot
}
\makeatother

\begin{document}

\begin{frontmatter}
\title{A data-driven and model-based accelerated Hamiltonian Monte Carlo method  for Bayesian elliptic inverse problems}

   \author[hku]{Sijing Li}
   \ead{lsj17@hku.hk}  
   \author[pku]{Cheng Zhang\corref{cor1}}
   \ead{chengzhang@math.pku.edu.cn}    
   \author[hku]{Zhiwen Zhang\corref{cor1}}
   \ead{zhangzw@hku.hk}
   \author[uci]{Hongkai Zhao}
   \ead{zhao@math.duke.edu}  
 
  \address[hku]{Department of Mathematics, The University of Hong Kong, Pokfulam Road, Hong Kong SAR, China.}
  \address[pku]{School of Mathematical Sciences and Center for Statistical Science, Peking University, Beijing, China.}
  \address[uci]{Department of Mathematics, Duke University, Durham, NC 27708, USA.}
  \cortext[cor1]{Corresponding author}
\begin{abstract} 
In this paper, we consider a Bayesian inverse problem modeled by elliptic partial differential equations (PDEs). Specifically, we propose a data-driven and model-based approach to accelerate the Hamiltonian Monte Carlo (HMC) method in solving large-scale Bayesian inverse problems. The key idea is to exploit (model-based) and construct (data-based) the intrinsic approximate low-dimensional structure of the underlying problem which consists of two components -- a training component that computes a set of data-driven basis to achieve significant dimension reduction in the solution space, and a fast solving component that computes the solution and its derivatives for a newly sampled elliptic PDE with the constructed data-driven basis. Hence we achieve an effective data and model-based approach for the Bayesian inverse problem and overcome the typical computational bottleneck of HMC -- repeated evaluation of the Hamiltonian involving the solution (and its derivatives) modeled by a complex system, a multiscale elliptic PDE in our case.  We present numerical examples to demonstrate the accuracy and efficiency of the proposed method.  
 


\end{abstract}
\end{frontmatter}

\section{Introduction} \label{sec:introduction}
\noindent
Inverse problems are ubiquitous in models used in science and engineering where problem-specific parameters or inputs need to be estimated from indirect and noisy observations. However, inverse problems are often nonlinear (even if the forward problems are linear) and ill-posed (or unstable) in that either the existence and uniqueness of the solutions may not be guaranteed or the dependence of the parameters on the data (and noise) may be sensitive. As a result, pointwise estimates may be erroneous and misleading and additional regularization is often required. On the other hand, the Bayesian approach to inverse problems \cite{kaipio:2005,dashti2011uncertainty,beskos2015sequential,lan2019adaptive} can provide another alternative.
In the Bayesian paradigm, the solution to the inverse problem is posited as the posterior distribution of the unknowns conditioned on observations where regularization is naturally imposed in the form of an appropriate prior distribution.
Bayesian inversion, therefore, provides a principled way of uncertainty quantification in the presence of data and noise.

As the posterior is generally intractable due to the complexity of the system, people often resort to computational approximation approaches such as Markov chain Monte Carlo (MCMC) methods.
In a typical MCMC method, samples from the posterior distribution are generated by updating the current states according to a proposing mechanism and a correction criterion designed to keep the posterior invariant. 
The efficiency of MCMC methods heavily depends on the design of proposing mechanism, i.e. its computation cost, acceptance probability, and mixing property.
For complicated and large systems in practice, it is important to strike an appropriate balance among these factors.
For example, simple MCMC algorithms (e.g., generating proposals based on random walk Metropolis), although easy and cheap to implement, usually have a low acceptance rate and mix poorly for complex and high dimensional problems since no information or structure of the underlying problem is utilized. 

In recent years, many advanced MCMC methods have been proposed to improve the sampling efficiency for high-dimensional problems \cite{duane87,neal2011mcmc}.
Based on an intelligent design of Hamiltonian dynamics, the HMC method uses gradient information of the underlying posterior distribution to make distant and less correlated proposals with high acceptance probabilities, greatly improving the mixing rate of the Markov chains.
On the other hand, the associated computation cost of those advanced MCMC methods can be a bottleneck that makes it difficult to scale up to complicated models and large data.
Note that the sampling procedure requires repetitive evaluations of the likelihood function and its derivatives and maybe other geometric and statistical quantities, e.g., Fisher information for Riemannian Hamiltonian Monte Carlo (RHMC) method \cite{girolami2011riemann}. 
To alleviate this issue, one popular attempt is to find a computationally cheap surrogate approximation to replace the original Hamiltonian \cite{zhang2017hamiltonian,zhang2017precomput,zhang2018vhmc,strathmann2015gradient,lan2016} in the sampling process. The key in designing an effective surrogate function is to capture the collective property of large datasets while removing redundancy. The overall computation efficiency is improved due to significant cost reduction in the Hamiltonian proposing process and insignificant loss in acceptance rate. Although these surrogate approximations can provide significant empirical performance improvement, they are usually obtained as blackbox approximations by fitting the training data where the mechanism and structure of the model that generates the data have largely remained unexplored and unexploited. 

For our Bayesian inverse problem, not only the unknown quantity is a random field that lives in high dimensions after discretization or can be a parametric model with many parameters, the solution to an elliptic PDE and its derivative are also involved in generating the data and evaluating the posterior and Hamiltonian in the HMC method.
These computational challenges make traditional MCMC methods extremely costly to use for Bayesian inverse problems.
In this work, we propose a data-driven  and model-based approach that can significantly reduce the computation cost of the HMC method  for Bayesian elliptic inverse problems.
The key idea is to exploit the intrinsic approximate low dimensional structure of elliptic differential operators and construct a data-driven basis as proposed in \cite{LiZhangZhao2019}.
First, a set of data-driven basis functions are constructed from training data, e.g., from real measurements or the initial burn-in stage of MCMC methods, to achieve significant dimension reduction in the solution space.
With the constructed basis, a newly sampled elliptic PDE can be solved efficiently.
Note that the derivatives (with respect to some parameters) of a solution to a linear PDE satisfies the same PDE (with different righthand sides) that can be computed efficiently as well.
Hence, this model-based and data-driven strategy can reduce the computation cost of the HMC sampling for our Bayesian inverse problem significantly. 



The rest of the paper is organized as follows. We first describe the forward model and the Bayesian inversion problem in Section~\ref{sec:modelproblem} and the HMC method for Bayesian inversion in Section~\ref{sec:HMCBayesianInversion}.  
Intrinsic low dimensional structure of the forward problem, model-based and data-driven dimension reduction, and approximation of the parameter-to-solution map is discussed in Section~\ref{sec:low-dim}.
The accelerated HMC (AHMC) method for Bayesian inverse problems is presented with implementation details in Section~\ref{sec:AHMC}. We present numerical experiments and results of AHMC and compare its performance to other state-of-the-art HMC methods in Section~\ref{sec:NumericalResults}.  Concluding remarks are made in Section~\ref{sec:conclusion}.
 
\section{Model problem} \label{sec:modelproblem} 
\subsection{Forward problem} \label{sec:forwardproblem} 
\noindent 
In this paper, we consider a classical inverse problem that involves inference of the diffusion coefficient in an elliptic PDE that is commonly used to model isothermal steady flow in porous media, hydrology and reservoir simulation, and many other applications. To be specific, we consider the following elliptic PDEs with random coefficients $a(\textbf{x},\omega)$, where one would like to infer, as the forward model,
\begin{align}
\mathcal{L}(\textbf{x},\omega) u(\textbf{x},\omega) 
\equiv -\nabla\cdot\big(a(\textbf{x},\omega)\nabla u(\textbf{x},\omega)\big) &= f(\textbf{x}), 
\quad \textbf{x}\in D, \quad \omega\in\Omega,   \label{MsStoEllip_ModelEq}\\
u(\textbf{x},\omega)&= 0, \quad \quad \textbf{x}\in \partial D, \label{MsStoEllip_ModelBC}
\end{align}
where $D \in \mathbb{R}^d$ is a bounded spatial domain, $\Omega$ is a sample space, and the source function $f(\textbf{x})\in L^2(D)$.  We assume $a(\textbf{x},\omega)$ in  \eqref{MsStoEllip_ModelEq} is  almost surely uniformly elliptic, namely, there exist $a_{\min}, a_{\max}>0$, such that
\begin{align}
P\big(\omega\in \Omega: a(\textbf{x}, \omega)\in [a_{\min},a_{\max}], \forall \textbf{x} \in D\big) = 1.
\label{asUniformlyElliptic1}
\end{align}
In general, we can assume the random coefficient $a(\textbf{x},\omega)$ is of some parametric form. For example, a commonly used affine form is the following,  
\begin{align}
a(\textbf{x},\omega) =  \bar{a}(\textbf{x}) + \sum_{m=1}^{r}a_m(\textbf{x})\xi_{m}(\omega),
\label{AffineParametrizeRandomCoefficient}
\end{align}
where $\xi_{m}(\omega)$, $m=1,...,r$ are random variables and $a_m(\textbf{x})$ are some spatial basis functions, e.g., finite element basis, polynomial basis, Fourier basis, radial basis, etc.  

Once a parametric form of the random coefficient $a(\textbf{x},\omega)$ is given, 
computing the solution $u(\textbf{x},\omega)$ to the problem \eqref{MsStoEllip_ModelEq}-\eqref{MsStoEllip_ModelBC} 
defines a map from the parameter domain $\vxi(\omega)=\big(\xi_1(\omega),\cdots,\xi_r(\omega)\big)^T \in \mathcal{W}\subset \mathbb{R}^r$ to the solution space 
\begin{align}
\vxi(\omega) \mapsto  u(\textbf{x},\omega) = u(\textbf{x},\vxi(\omega)) \in H_0^1(D),
\label{solution-map}
\end{align}
which is a Banach-space-valued function of the random input vector $\vxi(\omega)$. 

Many efficient numerical methods have been developed for solving elliptic PDEs with random coefficients; see \cite{Ghanem:91,Xiu:03,Zabaras:06,babuska:04, Babuska:07,Webster:08, graham2011quasi,abdulle2013multilevel, Grahamquasi:2015} and references therein. By solving the forward problem, one can quantify the uncertainty in the elliptic PDEs with randomness. However, when the elliptic PDEs involve multiscale features and/or high-dimensional random inputs, these problems become challenging due to high computational costs. 
In recent years, we have developed data-driven methods to solve multiscale elliptic PDEs with random coefficients  \eqref{MsStoEllip_ModelEq} based on intrinsic dimension reduction \cite{ZhangCiHouMMS:15,efendiev2015multilevel,chung2018cluster}. We also refer the intertested reader to \cite{Zabaras:13,abdulle2013multilevel,ZhangHouLiu:15,efendiev2015multilevel,chung2018cluster} for other methods to solve \eqref{MsStoEllip_ModelEq}.

\subsection{Bayesian inverse problems} \label{sec:inverseproblem}
\noindent
Let $\mathcal{W}$ be the space of admissible unknowns and $\mathcal{F}: \mathcal{W}\mapto \mathcal{U}$ be a forward map represents a mathematical model that assigns an output $u\in\mathcal{U}$ to an input $\vxi\in\mathcal{W}$.
In this paper, we focus on the elliptic PDE \eqref{MsStoEllip_ModelEq} where $\vxi$ is the parameters in the random coefficient $a(\textbf{x},\vxi)$ and $u$ is the solution to the PDE with the corresponding coefficient.
The inverse problem is to recover the unknown parameter $\vxi\in\mathcal{W}$ (and hence the coefficient $a(\textbf{x},\vxi)$) from some measurement of solution $u$ in the domain and at the boundary. Often in practice $u$ can only be recorded at finite discrete locations with noise which is the data denoted by $\textbf{y}\in \mathbb{R}^m$ related by 
\begin{equation}\label{eq:inverseprob}
\textbf{y} = \mathcal{G}(\vxi) + \bm{\eta}.
\end{equation}
Here the forward model $\mathcal{G}: \mathbb{R}^r\mapto \mathbb{R}^m$ is a composition of the forward map $\mathcal{F}$ and a discretized observation operator through which observable quantities (e.g., point-wise evaluation of the solution) are collected, and $\bm{\eta}\in\mathbb{R}^m$ is the measurement error (or the noise).

In the Bayesian formulation of the inverse problem \eqref{eq:inverseprob}, one treats the parameter $\vxi$ as a random variable (vector) with a prior distribution $p_{\vxi}(\vxi)$.
The noisy model, i.e., distribution of $\bm{\eta}$, gives the likelihood $p_{\textbf{y}|\vxi}(\textbf{y}|\vxi)$.
For simplicity and concreteness, in this paper we assume that $\bm{\eta}$ is a zero-mean Gaussian with diagonal covariance $\sigma^2\bm{I}_m$, so that
\begin{equation}\label{eq:likelihood}
p_{\textbf{y}|\vxi}(\textbf{y}|\vxi) \propto \exp\left(-\Phi(\vxi;\textbf{y})\right),\qquad \Phi(\vxi;\textbf{y}) := \frac{\|\textbf{y}-\mathcal{G}(\vxi)\|^2}{2\sigma^2}.
\end{equation}
The posterior distribution of $\vxi$ conditioned on the data $\textbf{y}$ then follows the Bayes' rule:
\begin{equation}
p_{\vxi|\textbf{y}}(\vxi|\textbf{y}) \propto p_{\textbf{y}|\vxi}(\textbf{y}|\vxi)\cdot p_{\vxi}(\vxi)
\end{equation}
and Bayesian inversion can be performed by estimating the posterior via, e.g., the HMC method and other MCMC methods.
 

In addition to the usual computational issues for MCMC type of methods, there is another challenge for the Bayesian elliptic inverse problem due to the complicated forward model \eqref{MsStoEllip_ModelEq}. Instead of a simple explicit probabilistic model that prescribes the likelihood of data given the parameter of interest, one needs to solve the elliptic PDE \eqref{MsStoEllip_ModelEq} for each coefficient corresponding to a new sample of the parameter $\vxi$ to compute the likelihood function \eqref{eq:likelihood}, which is the computation bottleneck for the Bayesian inversion. To address these challenges, we propose a data-driven and model-based accelerated HMC method that improves the convergence rate of the MCMC method and exploit the underlying forward model \eqref{MsStoEllip_ModelEq} using a data-driven approach proposed in \cite{LiZhangZhao2019}, which enables us to reduce the computational cost in solving the forward model problem and hence the overall sampling cost.

\section{The HMC method for Bayesian inversion} \label{sec:HMCBayesianInversion} 
\noindent
The HMC method is one of the state-of-the-art MCMC methods suitable for complex high dimensional target distributions with strong dependencies between parameters, which is the case for Bayesian inverse problems. Leveraging geometric information from the target distribution, the HMC method \cite{duane87,neal2011mcmc} extends the parameter space with auxiliary momentum variables $\bm{\zeta}$, and introduces a Hamiltonian dynamics system to propose samples of model parameters within the Metropolis framework, greatly enhancing the exploration efficiency in the parameter space compared to simple random walk proposals. More specifically, the HMC method generates proposals jointly for $\vxi$ and $\bm{\zeta}$ using the following system of differential equations
\begin{equation}\label{eq:HamiltonianDynamics}
\frac{d\vxi}{dt} = \frac{\partial H}{\partial \bm{\zeta}},\qquad
\frac{d\bm{\zeta}}{dt} = -\frac{\partial H}{\partial \vxi}.
\end{equation}
where the Hamiltonian function is defined as $H(\vxi, \bm{\zeta}) = U(\vxi) + K(\bm{\zeta})$.
Here in the Bayesian elliptic inverse problem, the potential energy $U$ is defined as $U(\vxi) = -\log p_{\textbf{y}|\vxi}(\textbf{y}|\vxi) - \log p_{\vxi}(\vxi)$, and the kinetic energy $K(\bm{\zeta}) = \frac12\bm{\zeta}^TM^{-1}\bm{\zeta}$ corresponds to the negative log-density of a zero-mean multivariate Gaussian distribution with covariance $M$ (also known as the mass matrix and is often set to be the identity).
As the analytical solution of the Hamiltonian dynamics \eqref{eq:HamiltonianDynamics} is usually unavailable, proposals in the HMC method are often made by numerical simulation via the leap-frog scheme. Speficically, given the sample $(\vxi^{(t)}, \bm{\zeta}^{(t)})$ at time $t$, we generate the sample at time $t+1$ by the following scheme
 \begin{equation}\label{eq:Leapfrog}
\begin{aligned}
\bm{\zeta}^{(t+\frac12)} &= \bm{\zeta}^{(t)} - \frac{\Delta t}{2}\nabla_{\vxi}U(\vxi^{(t)}),\\
\vxi^{(t+1)} &= \vxi^{(t)} + \Delta t \nabla_{\bm{\zeta}} K(\bm{\zeta}^{(t+\frac12)}),\\
\bm{\zeta}^{(t+1)} &= \bm{\zeta}^{(t+\frac12)} - \frac{\Delta t}{2}\nabla_{\vxi}U(\vxi^{(t+1)}),
\end{aligned}
\end{equation}
where $\Delta t$ is the step size. Starting from the current state $(\vxi, \bm{\zeta})$, where $\vxi$ is the current parameter and $\bm{\zeta}$ is resampled from the multivariate Gaussian distribution $\mathcal{N}(\bm{0}, M)$, the proposed state $(\vxi^\ast, \bm{\zeta}^\ast)$ at the end of a simulated trajectory of length $L$ is accepted with probability  
\begin{equation} 
p=\min\big(1, \exp[-H(\vxi^\ast, \bm{\zeta}^\ast) + H(\vxi, \bm{\zeta})]\big).
\label{acceptP}
\end{equation}	
From this point of view, the HMC method can be viewed as a Metropolis algorithm that samples from the joint distribution
\begin{equation}
p(\vxi,\bm{\zeta}) \propto \exp\left(-U(\vxi)-\frac12\bm{\zeta}^TM^{-1}\bm{\zeta}\right).
\end{equation}
The marginal distribution of $\vxi$ then follows the target posterior distribution since $\vxi$ and $\bm{\zeta}$ are separated (i.e., independent). Note that the Hamiltonian is preserved for analytical solutions of \eqref{eq:HamiltonianDynamics}, and the discretization error in \eqref{eq:Leapfrog} can be controlled by appropriate choice of the step size $\Delta t$, the HMC method is often able to generate distant, uncorrelated proposals with a high acceptance probability, allowing for efficient exploration of the parameter space.

For our Bayesian inverse problem, however, there is still a computational bottleneck we have to resolve, that is repetitive computation of solution $u(\textbf{x},\vxi)$ to the elliptic PDE \eqref{MsStoEllip_ModelEq} to evaluate the potential energy $U(\vxi) = -\log p_{\textbf{y}|\vxi}(\textbf{y}|\vxi) - \log p_{\vxi}(\vxi)$ in the Hamiltonian, and even more, the gradient with respect to the parameter $\nabla_{\vxi}U(\vxi)$ needs to be repetitively evaluated to simulate a trajectory for HMC proposals as in \eqref{eq:Leapfrog}. Note that the key to evaluation of $\nabla_{\vxi}U(\vxi)$ is the evaluation of derivatives $\frac{\partial u(\textbf{x},\vxi)}{\partial \xi_j}=u_{\xi_j}(\textbf{x},\vxi)$ of the solution to \eqref{MsStoEllip_ModelEq} satisfying
\begin{align}
 -\nabla\cdot\big(a(\textbf{x},\vxi)\nabla u_{\xi_j}(\textbf{x},\vxi)\big) &= \nabla \cdot \big(a_{\xi_j}(\textbf{x},\vxi) \nabla u(\textbf{x},\vxi)\big), 
\quad \textbf{x}\in D, \label{eq:gradient}\\
u_{\xi_j}(\textbf{x},\vxi)&= 0, \quad \quad \textbf{x}\in \partial D, \label{eq:bc}
\end{align}
which is the same elliptic PDE as \eqref{MsStoEllip_ModelEq} with a righthand side that depends on the solution to  \eqref{MsStoEllip_ModelEq} corresponding to the current sample of $\vxi$. This could easily become prohibitively expensive in practice since so many PDEs have to be solved for each sampling step.

In what follows, we describe how to approximate the low-dimensional structure of the solution space to the \eqref{MsStoEllip_ModelEq} with varying coefficients and righthand side and the data-driven approach proposed in \cite{LiZhangZhao2019} that can take advantage of the approximate low-dimensional structure of the forward model to accelerate the HMC method.

\section{Low dimensional structure of the forward model and approximation of the parameter-to-solution map}\label{sec:low-dim}
\noindent
For the Bayesian inverse problem, we are facing the challenge to solve \eqref{MsStoEllip_ModelEq} with different coefficients and different righthand sides \eqref{eq:gradient} repetitively in the sampling process.

\subsection{Low dimensional structure with respect to varying coefficients}
\noindent
With the uniform ellipticity assumption of $a(\textbf{x},\vxi(\omega))$ and its smooth dependence on the parameter $\vxi$, 
the solution $u(\textbf{x},\vxi)$ also depends smoothly on the parameters, which can be approximated via a polynomial expansion in $\vxi$ of the form 
\begin{align}
\sum_{\gvec{\alpha}\in\mathcal{J}_{r}}u_{\gvec{\alpha}}(\textbf{x})\vxi^{\gvec{\alpha}}(\omega),
\label{taylorseries} 
\end{align}
where $\gvec{\alpha}=(\alpha_1, \alpha_2, \cdots, \alpha_{r})$ is a multi-index, $\mathcal{J}_{r}=\{\gvec{\alpha}\,|\, \alpha_i \ge 0, \alpha_i \in \mathbb{N}, 1\leq i \leq r\}$ is a multi-index set of countable cardinality, and  $\vxi^{\gvec{\alpha}}(\omega)=\prod_{1\leq i \leq r}\xi_{i}^{\alpha_i}(\omega)$ is a multivariate polynomial. 
 
In particular, if uniform ellipticity assumption of $a(\textbf{x},\vxi)$ has a holomorphic extension to an open set in complex domain that contains the real domain for $\vxi$, explicit estimates for the coefficients $u_{\gvec{\alpha}}$ can be established similar to those estimates for the polynomial approximation for an analytic function. From the estimates, the following result for the best $n$-term approximation can be proved (see \cite{cohen2015approximation} for details).

\begin{proposition}\label{BestnTermAppx} 
	Consider a parametric problem of the form \eqref{MsStoEllip_ModelEq}-\eqref{MsStoEllip_ModelBC} with a random coefficient \eqref{AffineParametrizeRandomCoefficient}. Both the Taylor series and Legendre series of the form \eqref{taylorseries} converges to $u(x,\vxi(\omega))$ in $H_0^1(D)$ for all $\vxi(\omega)\in \mathcal{W}$. Moreover, for any set $\mset{J}_{r}^{n}$ of indices corresponding to the $n$ largest of $||u_{\gvec{\alpha}}(\cdot)||_{H_0^1(D)}$, we have   	
	\begin{align}
	\sup_{\vxi(\omega)\in\mathcal{U}}\big|\big|u(\cdot,\vxi(\omega))-\sum_{\gvec{\alpha}\in\mset{J}_{r}^{n}}u_{\gvec{\alpha}}(\cdot)\vxi^{\gvec{\alpha}}(\omega)\big|\big|_{H_0^1(D)} \leq C \exp(-c n^{1/r}),\label{taylorseries-n-bestterms} 
	\end{align}
	where $\mathcal{J}_{r}^{n}$ is a subset of $\mathcal{J}_{r}$ with cardinality $\#\mathcal{J}_{r}^{n}=n$, $C$ and $c$ are positive and depend on $r$. 
\end{proposition}

Proposition \ref{BestnTermAppx} shows that there exists a linear subspace with dimension at most $O(n\sim(\frac{\log C}{c}+\frac{|\log\epsilon| }{c})^r)$, e.g., spanned by $u_{\gvec{\alpha}}(x)$, $\gvec{\alpha}\in\mathcal{J}_{r}^{n}$, that can approximate the solution of \eqref{MsStoEllip_ModelEq}-\eqref{MsStoEllip_ModelBC} with random coefficient within $\epsilon$ error.


The result in proposition \ref{BestnTermAppx} reveals the existence of approximate low dimensional structures in the solution space of \eqref{MsStoEllip_ModelEq}-\eqref{MsStoEllip_ModelBC}. However, this approximation is obtained by mathematical analysis, which cannot be directly implemented via a computational algorithm. In \cite{LiZhangZhao2019}, a data-driven approach was proposed to construct problem dependent basis functions that can approximate the solution space of \eqref{MsStoEllip_ModelEq}-\eqref{MsStoEllip_ModelBC} effectively. 

\begin{remark}
  When the coefficient $a(\textbf{x},\omega)$ is a nonlinear function of a finite number of random variables, one can apply the empirical interpolation method (EIM) \cite{PateraMaday:2004} to approximately convert $a(\textbf{x},\omega)$ into an affine form. Thus, low dimensional structures still exist in the solution space. In addition, we refer the reader to \cite{hoang2014n,bachmayr2017sparse} for the results of the best $n$-term polynomial approximation of elliptic PDEs with lognormal coefficients.
\end{remark} 


 \subsection{Low dimensional structure with respect to varying sources}
 \noindent
 Consider the following elliptic PDE with a fixed coefficient $a(\textbf{x}) \in L^{\infty}(D)$ and a random source $f(\textbf{x},\omega)$,
 \begin{align}
\label{eq:source}
 -\nabla\cdot\big(a(\textbf{x})\nabla u(\textbf{x},\omega)\big) &= f(\textbf{x},\omega), 
\quad \textbf{x}\in D, \\
u(\textbf{x},\omega)&= 0, \quad \quad \textbf{x}\in \partial D, 
\end{align}
with $a_{\min}, a_{\max}>0$, such that $a_{\min}<a(\textbf{x})<a_{\max}$ for all $\textbf{x} \in D$. The contrast ratio $\kappa_a=\frac{a_{\max}}{a_{\min}}$ is an important factor in the stability and convergence analysis. Let $G(\textbf{x},\textbf{y})$ be the Green's function for the elliptic PDE \eqref{eq:source} satisfying 
\begin{equation}
 -\nabla\cdot\big(a(\textbf{x})\nabla G(\textbf{x},\textbf{y})\big) = \delta(\cdot,\textbf{y}), \quad \text{in}~ D,  \quad  G(\cdot,\textbf{y})=0, \quad \text{on} ~ \partial D,
\end{equation}
where $\delta(\cdot,\textbf{y})$ is the Dirac delta function denoting an impulse source point at $\textbf{y}\in D$. The solution to \eqref{eq:source} can be represented as
\begin{equation}
u(\textbf{x},\omega)=\int_{D}G(\textbf{x},\textbf{y})f(\textbf{y},\omega)d\textbf{y}. \label{solutionrepresentationGreen}
\end{equation}

It was shown in \cite{BebendorfHackbusch:2003} that the Green's function for an elliptic operator is highly separable.
\begin{proposition}[Theorem 2.8 of \cite{BebendorfHackbusch:2003}]\label{GreenFuncSepaApp}
	
	Let $D_1, D_2 \subset D$ be two subdomains and $D_1$ be convex. Assume that there exists $\rho>0$ such that 
	\begin{align}
	0 < \text{ \normalfont diam} (D_1) \leq \rho\text{ \normalfont dist} (D_1, D_2). 
	\label{AdmissiblePairs}
	\end{align}
	Then, for any $\epsilon \in (0,1)$ there is a separable approximation
	\begin{align}
	G_k(\textbf{x},\textbf{y}) = \sum_{i=1}^k u_i(\textbf{x}) v_i(\textbf{y}) \quad \text{with } k \leq  
	c^d(\kappa_a, \rho) |\log \epsilon|^{d+1},
	\label{GreenFuncSepaApp1}
	\end{align}
	so that for all $\textbf{y}\in D_2$
	\begin{align}
	\| G( \cdot,\textbf{y}) - G_k(\cdot,\textbf{y}) \|_{L^2 (D_1)} \leq \epsilon \| G(\cdot,\textbf{y}) \| _{L^2(\hat{D}_1)},
	\end{align}
	where  $\hat{D}_1 := \{ \textbf{x} \in D :  2\rho~\text{\normalfont dist} (\textbf{x}, D_1) \leq \text{\normalfont diam} (D_1)\}$.
\end{proposition}

The above result shows that there exists a low dimensional linear subspace, e.g., spanned by $u_i(\cdot)$, that can approximate the family of functions $G(\cdot,\textbf{y})$ well in $L^2(D_1)$ uniformly with respect to $\textbf{y}\in D_2$. Moreover, if $\mathrm{supp}(f(\textbf{x},\theta))\subset D_2$,  one can approximate the family of solutions $u(\textbf{x},\omega)$ to \eqref{eq:source} by the same space well in $L^2(D_1)$ uniformly. Indeed, let
 \begin{equation}
u^{\epsilon}(\textbf{x},\omega)=\int_{D_2} G_k(\textbf{x},\textbf{y})f(\textbf{y},\omega) d\textbf{y}=\sum_{i=1}^k u_i(\textbf{x})\int_{D_2} v_i(\textbf{y}) f(\textbf{y},\omega) d\textbf{y}.
\end{equation}
We have
\begin{equation}
\begin{array}{l}
\|u(\cdot,\omega)-u^{\epsilon}(\cdot,\omega)\|^2_{L^2(D_1)}=\int_{D_1} \left[\int_{D_2} (G(\textbf{x},\textbf{y})-G_k(\textbf{x},\textbf{y}))f(\textbf{y},\omega) d\textbf{y}\right]^2 d\textbf{x} 
\\ \\
\le \|f\|_{L^2(D_2)}^2 \int_{D_2}\| G(\cdot,\textbf{y}) - G_k(\cdot,\textbf{y}) \|^2_{L^2 (D_1)} d\textbf{y}\le C(D_1, D_2, \kappa_a, d)\epsilon^2\|f\|_{L^2(D_2)}^2,
\end{array}
\end{equation}
since $\| G(\cdot,\textbf{y}) \| _{L^2(\hat{D}_1)}$ is bounded uniformly with respect to $\textbf{y}\in D_2$ by a positive constant that depends on $D_1, D_2, \kappa_a, d$ due to the uniform ellipticity. Note that the low dimensional structure does not need any regularity assumption in $a(\textbf{x})$. Moreover, dependence of the source on randomness can be arbitrary in terms of dimensionality and regularity.
\begin{remark}
Although, the proof of high separability of the Green's function requires $\textbf{x}\in D_1, \textbf{y}\in D_2$  for two disjoint $D_1$ and $D_2$ due to the singularity of the Green's function at $\textbf{x}=\textbf{y}$,  the above approximation of the solution $u$ in a domain disjoint with the support of $f$ also works for $u$ in the whole domain even when $f$ is a globally supported smooth function as shown in our numerical results in \cite{LiZhangZhao2019}. 
\end{remark}


\subsection{Data-driven basis for dimension reduction} \label{sec:ConstructDDbais}
\noindent
Since there exist low-dimensional structures in the solution space of elliptic PDEs with random coefficients and sources, we use problem-specific and data-driven basis to achieve a significant dimension reduction in solving the elliptic PDEs \eqref{MsStoEllip_ModelEq}. Our method consists of a training process and a solving process. In the training process, we extract the low-dimensional structure of the solution space and construct a set of data-driven basis functions from training data or real measurements, e.g., a set of solution samples $\{u(\textbf{x},\omega_i)\}_{i=1}^{N}$ can be obtained from measurements or generated by solving  the elliptic PDE \eqref{MsStoEllip_ModelEq}-\eqref{MsStoEllip_ModelBC}, e.g., with coefficient samples $\{a(\textbf{x},\omega_i)\}_{i=1}^{N}$ during the burning stage of the HMC method.  

Let $V_{snap}=\{u(\textbf{x},\omega_1),...,u(\textbf{x},\omega_N)\}$ denote the solution samples. 
We use the POD method \cite{HolmesLumleyPOD:1993,Sirovich:1987,Willcox2015PODsurvey}, or a.k.a the PCA method, to find the optimal subspace and its orthonormal basis functions to approximate $V_{snap}$ to a certain accuracy. Specifically, we define the correlation matrix 
$\Sigma=(\sigma_{ij})\in \mathbb{R}^{N\times N}$ with  $\sigma_{ij}=<u(\cdot,\omega_i), u(\cdot,\omega_j)>_{D}$, $i,j= 1, \ldots, N$, where $<\cdot, \cdot>_{D}$ denotes the standard inner product on $L^2(D)$. Let the eigenvalues of the correlation matrix be $\lambda_1\ge  \lambda_2 \ge \ldots \ge \ldots \ge \lambda_N \ge 0$ and the corresponding eigenfunctions be $\phi_{1}(x)$, $\phi_{2}(x), \ldots, \phi_N(x)$, which will be referred to as data-driven basis functions.
\begin{proposition}\label{POD_proposition}
	 The space spanned by the leading $K$ data-driven basis functions has the following approximation property to $V_{snap}$. 
	\begin{align}
	\frac{\sum_{i=1}^{N}\Big|\Big|u(x,\omega_{i})- \sum_{j=1}^{K}<u(\cdot,\omega_{i}),
		\phi_j(\cdot)>_{D}\phi_j(x)\Big|\Big|_{L^2(D)}^{2} }{\sum_{i=1}^{N}\Big|\Big|u(x,\omega_{i})\Big|\Big|_{L^2(D)}^{2}}=\frac{\sum_{s=K+1}^{N}  \lambda_s}{\sum_{s=1}^{N}  \lambda_s}.
	\label{Prop_PODError}
	\end{align}
\end{proposition} 
First, we expect a fast decay in the eigenvalues $\lambda_s$ so that a small set of data-driven basis functions ($K\ll N$) will be enough to approximate the solution samples well in the root mean square sense. Secondly, based on the existence of low-dimensional structure, we expect that the data-driven basis functions, $\phi_{1}(\textbf{x})$, $\phi_{2}(\textbf{x}), \ldots, \phi_{K}(\textbf{x})$,  can approximate the solution $u(\textbf{x},\omega)$ well by 
$u(\textbf{x},\omega) \approx \sum_{j=1}^{K}c_{j}(\omega)\phi_{j}(\textbf{x})$ almost surely for $\omega \in \Omega$.

Determining a set of good solution samples is important for the construction of the data-driven basis functions. In general, this issue is very challenging especially when the dimension of the random coefficient is high. Under
certain assumptions on the random coefficient, we obtained some criteria on how to choose the coefficient samples in order to obtain a set of accurate data-driven basis functions; see Section 3.4 of  \cite{LiZhangZhao2019}.

The computational costs of constructing the data-driven basis functions consist of two parts, if data are generated by simulation: (1) compute solution samples $\{u(\textbf{x},\omega_i)\}_{i=1}^{N}$; and (2) compute the data-driven basis by the POD method. This is common nature for many model reduction methods. Effective samples of solutions (see Section 3.4 of  \cite{LiZhangZhao2019}) and the use of randomized algorithms  \cite{Tropp:11} for the singular value decomposition (SVD) (utilizing the low-rank structure) help to reduce the offline computation cost. 


Equipped with the data driven basis $\phi_{j}(\textbf{x})$, $j=1,...,K$, we  can solve the problem \eqref{MsStoEllip_ModelEq}-\eqref{MsStoEllip_ModelBC} on the domain $D$ by the standard Galerkin formulation for  new realizations of $a(\textbf{x},\omega)$. Specifically, given a new realization of the coefficient $a(\textbf{x},\omega)$, we approximate the corresponding solution $u(\textbf{x},\omega)$ as 
\begin{align}
	u(\textbf{x},\omega) \approx \sum_{j=1}^{K}c_{j}(\omega)\phi_{j}(\textbf{x}), \quad \text{a.s. }
	\omega \in \Omega,  \label{RB_expansion2}
\end{align} 
and use the Galerkin projection to determine the coefficients $c_{j}(\omega)$, $j=1,...,K$.
We substitute the approximation \eqref{RB_expansion2} into Eq.\eqref{MsStoEllip_ModelEq}, multiply both side by $\phi_{l}(\textbf{x})$, $l=1,...,K$,  take integration over the domain $D$, and obtain a coupled linear system as folows 
\begin{align}
	\sum_{j=1}^K \int_{D}a(\textbf{x},\omega)c_{j}(\omega)\nabla\phi_{j}(\textbf{x})\cdot\nabla\phi_{l}(\textbf{x})d\textbf{x}  = \int_{D}f(\textbf{x})\phi_{l}(\textbf{x})d\textbf{x}, 
	\quad l=1,...,K.  \label{GalerkinSystem}
\end{align}
The computational cost of solving the linear system \eqref{GalerkinSystem} is small compared to using a Galerkin method, such as the finite element method, directly for $u(\textbf{x},\omega)$ because $K$ is much smaller than the degree of freedom needed to discretize $u(\textbf{x},\omega)$ in the whole domain.   

Note that if $a(\textbf{x},\omega)$ has the affine form \eqref{AffineParametrizeRandomCoefficient}, 
we first compute the terms that do not depend on randomness, including 
$\int_{D}\bar{a}(\textbf{x})\nabla\phi_{j}(\textbf{x})\cdot\nabla\phi_{l}(\textbf{x})d\textbf{x}$,  $\int_{D}a_{m}(\textbf{x})\nabla\phi_{j}(\textbf{x})\cdot\nabla\phi_{l}(\textbf{x})d\textbf{x}$ and
$\int_{D}f(\textbf{x})\phi_{j}(\textbf{x})d\textbf{x}$, $j,l=1,...,K$. Then, we save them in the offline stage.  This leads to considerable savings in assembling the stiffness matrix for each new realization of the coefficient $a(\textbf{x},\omega)$ in the online stage.

\subsection{The parameter-to-solution map}\label{sec:parametrized}
\noindent
To solve the Bayesian inverse problem modeled by the elliptic PDE \eqref{MsStoEllip_ModelEq}, we need to compute $c_{j}(\omega)$ by solving the linear equation system \eqref{GalerkinSystem} for many realizations of $a(\textbf{x},\omega)$. Although the data-driven basis functions provide considerable saving over standard finite element basis functions in solving \eqref{MsStoEllip_ModelEq}, it still requires a certain amount of computational cost in solving the linear equation system \eqref{GalerkinSystem} in the HMC methods. To further reduce the computational cost in the HMC method, we construct parameter-to-solution maps based on the training solution data and the data-driven basis functions.

According to our assumption, $a(\textbf{x},\omega)$ is parameterized by $r$ independent random variables,  i.e., 
$a(\textbf{x},\omega) =  a(\textbf{x},\xi_{1}(\omega),...,\xi_{r}(\omega))$. 
Thus, the solution can be represented as a functional of these random variables as well, i.e., $u(\textbf{x},\omega) = u(\textbf{x},\xi_{1}(\omega),...,\xi_{r}(\omega))$.
Let $\vxi(\omega)=[\xi_1(\omega),\cdots,\xi_r(\omega)]^T$ denote the 
random input vector and $\textbf{c}(\omega)=[c_{1}(\omega),\cdots,c_{K}(\omega)]^T$ denote the vector of solution coefficients in \eqref{RB_expansion2}. Now,  the problem can be viewed as constructing 
a parameter-to-solution map from $\vxi(\omega)$ to $\textbf{c}(\omega)$, denoted by $\textbf{F}:\vxi(\omega)\mapsto \textbf{c}(\omega)$, which is nonlinear. We approximate this nonlinear map through the given solution or measurement data. Given a set of solution samples $\{u(\textbf{x},\omega_i)\}_{i=1}^{N}$  corresponding to $\{\vxi(\omega_i)\}_{i=1}^{N}$, e.g., by solving \eqref{MsStoEllip_ModelEq}-\eqref{MsStoEllip_ModelBC} with $a(\textbf{x},\xi_{1}(\omega_i),...,\xi_{r}(\omega_i))$,
from which the set of data driven basis $\phi_{j}(\textbf{x}), j=1,...,K$ is obtained by using POD method as described above,  we can easily compute the projection coefficients $\big\{\textbf{c}(\omega_i)\big\}_{i=1}^{N}$ of $u(\textbf{x},\omega_i)$ on $\phi_{j}(\textbf{x})$, $j=1,...,K$, i.e., $c_j(\omega_i)=<u(\textbf{x},\omega_i), \phi_{j}(\textbf{x})>_{D}$.  From the data set,
$F(\vxi(\omega_i))= \textbf{c}(\omega_i)$, $i=1,...,N$, we construct the map $\textbf{F}$. Note the significant dimension reduction by reducing the map $\vxi(\omega)\mapsto u(\textbf{x},\omega)$ to the map  $\vxi(\omega)\mapsto \textbf{c}(\omega)$.
We provide several ways to construct $\textbf{F}$, depending on the dimension of the random input vector. More implementation details can be found in \cite{LiZhangZhao2019}.

When the dimension of the random input $r$ is small or moderate, one can use interpolation. In particular, if the solution samples correspond to $\vxi$ located on a uniform or sparse grid, standard polynomial interpolation can be used to approximate the coefficient $c_j$ at a new point of $\vxi$. If the solution samples correspond to $\vxi$ at scattered points or the dimension of the random input $r$ is moderate or high, one can first find a few nearest neighbors to the new point efficiently using the $k$-$d$ tree algorithm \cite{wald2006building} and then use moving least square approximation centered at the new point to approximate the mapped value.

When the dimension of the random input $r$ is high, the interpolation approach becomes expensive and less accurate. Due to the dimension reduction by the data-driven basis functions, one can train a neural network with a small output dimension to approximate the parameter-to-solution map $\textbf{F}$. Numerical results in 
\cite{LiZhangZhao2019} show that this approach works well. 
We will adopt the neural network approach to approximate the parameter-to-solution maps for both the solution and its derivatives in this work. 
 
In the HMC method, one can compute the solution $u(\textbf{x},\omega)$ using the constructed map $\textbf{F}$. 
For example, given a new sample of $a(\textbf{x},\xi_{1}(\omega_i),...,\xi_{r}(\omega_i))$, we plug $\vxi(\omega)$ into the constructed map $\textbf{F}$ to approximate $\textbf{c}(\omega)=\textbf{F}(\vxi(\omega))$, which are the projection coefficients of the solution on the data-driven basis. So we can quickly obtain the new solution $u(\textbf{x},\omega)$ using Eq.\eqref{RB_expansion2}, where the computational time is negligible. 
Similarly, we can construct data-driven basis functions and approximate the parameter-to-solution maps for computing the partial derivatives of the solution. Once we obtain the numerical solutions and their derivatives, we can use them as a proposal in the HMC method. Numerical experiments show that our new method achieves significant savings in computing a new proposed sample over the standard HMC method.

\section{The accelerated HMC method and implementation} \label{sec:AHMC}
\noindent
In this section, we present the data-driven and model-based accelerated HMC method for solving Bayesian elliptic inverse problems with implementation details.

In the burn-in stage, we run the standard HMC method, i.e., solving the forward elliptic problem \eqref{MsStoEllip_ModelEq} for $u$ and solving \eqref{eq:gradient} for $u_{\xi_i}$ for the numerical evaluation of  Hamiltonian dynamics in \eqref{eq:Leapfrog} using standard finite element method. The samples of solution and its derivatives computed during the burn-in stage are collected and used to construct data-driven basis for dimension reductions using POD as described in Section \ref{sec:ConstructDDbais}. In particular, a set of basis is computed for $u$ and each $u_{\xi_i}$. Then we use the collected samples of solution and its derivatives to train two neural networks using the Adam optimization method (see \cite{kingma2014adam}) to approximate the paraemter-to-solution map described in Section \ref{sec:parametrized}. Although $u$ and $u_{\xi_j}$ satisfy the same elliptic PDE, $u$ has a fixed righthand source and $u_{\xi_i}$ has a varying righthand source. We find that it is more efficient and accurate to construct two separate neural networks to approximate the parameter-to-solution map, one for $u$ and one for all $u_{\xi_j}$. The neural network that approximates the parameter-to-solution map
for $u$ has a first layer that is a fully connected affine transform $\textbf{h}_{1}=\textbf{W}_{1}\vxi + \textbf{b}_{1}$. The following hidden layers are residual connections $\textbf{h}_{l}=\tanh(\textbf{W}_{l}\textbf{h}_{l-1} + \textbf{b}_{l}) + \textbf{h}_{l-1}$ (see \cite{he2016deep}). The output layer is another affine transform with output $\textbf{c}(\vxi)=(c_1(\vxi), c_2(\vxi), \ldots, c_K(\vxi))^T$ and the error to minimize is
$\sum_{j=1}^{N}\sum_{k=1}^{K}| c_k(\vxi_j)-\bar{c}_k(\vxi_j)|^2$, where $\bar{c}_k(\vxi_j), k=1, 2, \ldots, K$ is the projected coefficients from $j$-th data $
u(\textbf{x}, \vxi_j), j=1, 2, \ldots, N$, collected during the burn-in stage. 
The neural network that approximates the parameter-to-solution map for all $u_{\xi_i}, i=1, 2, \ldots, r$, where $r$ is the dimension of the parameter space, has a similar network structure as above with an output of $(\textbf{c}^1(\vxi),  \textbf{c}^2(\vxi), \dots, \textbf{c}^r(\vxi))$ and the error to minimize is $\sum_{i=1}^r\sum_{j=1}^{N}\sum_{k=1}^{K_i} |c^i_k(\vxi_j)-\bar{c}^i_k(\vxi_j)|^2$, where $\bar{c}^i_k(\vxi_j), k=1, 2, \ldots, K_i$ is the projected coefficients computed from $j$-th data $
u_{\xi_i}(\textbf{x}, \vxi_j), j=1, 2, \ldots, N$ collected during the burn-in stage.

Once the parameter-to-solution map is trained, we can evaluate the potential energy $U(\vxi)$ and its gradient $\nabla_{\vxi}U(\vxi)$ efficiently and hence significantly accelerate the HMC method to get the posterior samples by evolving the Markov chain as described in Section \ref{sec:HMCBayesianInversion} for Bayesian inverse problems. We list the implementation steps of the accelerated HMC algorithm in Algorithm \ref{Algorithm-acceleratedHMC}. 

\begin{algorithm}[h]   
	\caption{The accelerated HMC algorithm.} \label{Algorithm-acceleratedHMC}  
	\begin{algorithmic}[1]  
		\State Input: the prior distribution for $a(\textbf{x},\vxi(\omega))$. 
		\State Collect samples of solution and its partial derivatives, i.e. $\{\vxi_j, u(\vxi_j), \frac{\partial u(\vxi_j)}{\partial \xi_1}, \cdots,  \frac{\partial u(\vxi_j)}{\partial \xi_r}\}_{j=1}^{N}$ during the burn-in stage.
		\State Extract basis functions $\{\phi_j(\textbf{x})\}_{j=1}^{K}$ for the solution and basis functions $\{\phi^i_j(\textbf{x})\}_{j=1}^{K_i}$ for the partial derivatives of the solution, $i=1,\cdots, r$, from the collected data using the POD method.
		\State Get training data $\{\vxi_j, \textbf{c}(\vxi_j), \textbf{c}^1(\vxi_j), \cdots, \textbf{c}^r(\vxi_j)\}_{j=1}^{N}$ by projecting the samples of solution and its partial derivatives onto the corresponding basis.
		\State Train a neural network fitting the data pair $\{\vxi, \textbf{c}(\vxi)\}$, and train another neural networks fitting the data pair $\{\vxi, \textbf{c}^1(\vxi), \cdots, \textbf{c}^r(\vxi)\}$ to approximate the parameter-to-solution maps.
		\State Generate samples from posterior distribution via the data-driven accelerated HMC algorithm. 
		\begin{enumerate}
			\item[(1)] at the current position $\vxi$, sample a new momentum $\bm{\zeta}\sim\mathcal{N}(0,M)$ to get a starting point $(\vxi,\bm{\zeta})$;
			\item[(2)] apply the leap-frog scheme \eqref{eq:Leapfrog} to compute the Hamiltonian dynamic, with data-driven gradient for the potential via the learned parameter-to-solution maps in step 5;
			\item[(3)] accept the proposed sample $(\vxi^\ast,\bm{\zeta}^\ast)$ at the end of the trajectory with probability \eqref{acceptP}, where $H$ is computed using the FEM reference solution.
		\end{enumerate} 
		\State   Output: samples of $\{\vxi\}$  that converge to the posterior distribution.
	\end{algorithmic}  
\end{algorithm}  

The extra computation cost for our proposed accelerated HMC is the construction of data-driven basis and training of the two neural networks to approximate the parameter-to-solution map. Due to the intrinsic approximate low dimensional structure of the solution (and its derivatives) of the forward elliptic model, the singular values of the data covariance matrix decays very fast. 
So one only needs to approximate the space spanned by a few leading singular vectors which can be computed efficiently using randomized SVD algorithms as described in Section \ref{sec:ConstructDDbais} (and more details in \cite{LiZhangZhao2019}). After significant model based dimension reduction, a rather shallow and small neural network with simple structure and low dimension input and output is needed to approximate the parameter-to-solution map well in practice. Hence, evaluation of the constructed parameter-to-solution map vs a full computation of the forward elliptic PDE significantly reduces the computation cost in each leap-frog step of the HMC dynamic after the burn-in stage. Moreover, our data and model based dimension reduction captures the intrinsic low dimension structure of the underlying problem with a data-driven basis and accuracy control (through the POD) that strikes a good balance between computation efficiency (by dimension reduction and parameter-to-solution approximation) and exploration efficiency (by proposing well decorrelated samples with high acceptance rate), as demonstrated by numerical experiments in the next section.

\section{Numerical Experiments and Results}\label{sec:NumericalResults}
\noindent
In this section, we use numerical experiments to demonstrate the accuracy and efficiency of accelerated HMC method for Bayesian inverse problems, with comparison to other state-of-the-art methods, including the standard HMC method, and random network surrogate method \cite{zhang2017hamiltonian}.  The Python codes are published on GitHub.\footnote{https://github.com/LSijing/Bayesian-pde-inverse-problem.}

We consider the elliptic inverse problem 
\begin{align}\label{forwardex1}
-\nabla\cdot\big(a(\textbf{x},\omega)\nabla u(\textbf{x},\omega)\big) &= 0, \quad \textbf{x}=(x_1,x_2)\in [0,1]\times[0,1]
\end{align}
with mixed boundary condition,
\begin{align}\label{forwardBC}
\frac{\partial u(\textbf{x},\omega)}{\partial \bf{n}} |_{x_1=0, x_1=1} = 0, 
\quad u(\textbf{x},\omega)|_{ x_2=0}  = x_1, 
\quad u(\textbf{x},\omega)|_{x_2=1} = 1-x_1.
\end{align}
 
\subsection{A log-normal coefficient with isotropic heterogeneity}\label{sec:numericalex1}
\noindent
In the first example, a Gaussian prior with zero mean and covariance function
\begin{align}
c(\textbf{x},\textbf{x}^{\prime}) = \sigma_{a}^2\exp\Big(-\frac{||\textbf{x}-\textbf{x}^{\prime}||_2^2}{2l^2}\Big)
\label{CovarianceFunction1}
\end{align}
is assumed on $\log(a(\textbf{x},\omega))$, where $\textbf{x}$, $\textbf{x}^{\prime}$ are any two points on $[0,1]\times[0,1]$, and the parameters $\sigma_{a}^2$ and $l$ denote the variance and the correlation length, respectively. The diffusion coefficient is approximated via a truncated Karhunen-Lo\`eve (KL) expansion
\begin{align}
\log(a(\textbf{x},\vxi)) = \sum_{i=1}^{r} \xi_i\sqrt{\lambda_i}v_i(\textbf{x}),
\label{KLE1}
\end{align}
by $r$ i.i.d. Gaussian random variables $\xi_i$, where $\vxi=(\xi_1,...,\xi_r)$, $\lambda_i$ and $v_i(\textbf{x})$, $i=1,2,\cdots,r$ are eigenvalues and eigenfunctions of the prior covariance function \eqref{CovarianceFunction1}. In this experiment, we test the performance of the accelerated HMC algorithm~\ref{Algorithm-acceleratedHMC} for different input random dimensions, $r=25, 30, 35$. 

Suppose the observation $\textbf{y}=\big(y_1,y_2,\cdots, y_m\big)$ is obtained by adding independent Gaussian noise to the exact solutions at some measurable locations
\begin{align}
	y_j = u(x_j,\vxi) + \eta_j, \eta_j\sim \mathcal{N}(0,\sigma^2), \quad j=1,2,\cdots,m. 
	\label{observation1}
\end{align} 
Our goal is to infer $\vxi$ and hence $a(\textbf{x},\vxi)$ based on the  observation data $\textbf{y}$. In the Bayesian framework, the posterior on $\vxi$ is
\begin{align}
	p_{\vxi|\textbf{y}}(\vxi|\textbf{y}) &\propto p_{\textbf{y}|\vxi}(\textbf{y}|\vxi)\cdot p_{\vxi}(\vxi)\nonumber\\
	&\propto \exp(-\frac{1}{2\sigma^2}|\textbf{y}-u(\textbf{x},\vxi)|^2)\exp(-\frac{1}{2}\vxi^T\vxi),
	\label{posterior1}
\end{align}
which is the target distribution. Thanks to an efficient approximation to the  parameter-to-solution map, 
 the computation cost of $u(\textbf{x}, \vxi)$ and $\nabla_{\vxi} u(\textbf{x},\vxi)$ is significantly reduced and hence the computation of likelihood~\eqref{posterior1} and Hamiltonian dynamics~\eqref{eq:Leapfrog} in accelerated HMC is very fast. Moreover, as demonstrated later on, the acceptance rate and exploration efficiency do not compromise much as a consequence of model and data based dimension reduction. So the overall performance of the HMC method is enhanced significantly. 
 
To generate the training data, the discretization is done on a uniform grid with $31\times 31$ points through triangle finite element basis functions. Suppose the measurements are placed on $11\times 11$ grids of the numerical solution $u(\textbf{x},\cdot)$, i.e. $m=121$ in \eqref{observation1}. We choose $\sigma=0.1$ be the noise in the observation data \eqref{observation1}, and $\sigma_{a}=0.5$, $l=0.2$ be parameters in the prior covariance function \eqref{CovarianceFunction1}.

The burn-in stage consists of 10000 steps of standard HMC, of which 9000 accepted samples of solutions $u$ and its derivatives $\frac{\partial u}{\partial \xi_i}, i=1, 2, \ldots, r$ are collected during the burn-in stage. These collected data are first used to construct a set of data-driven basis using POD for dimension reduction. In our previous study \cite{LiZhangZhao2019}, we found that a larger number of basis functions are needed to approximate the derivatives of solution than those needed to approximate the solution.  Specifically, we construct $K=20$ basis $\phi_1(\textbf{x}_j), \phi_2(\textbf{x}_j), \ldots, \phi_K(\textbf{x}_j), j=1, 2, \ldots, m$ for the approximation of $u$ and $K_i=40$ basis $\phi^i_1(\textbf{x}_j), \phi^i_2(\textbf{x}_j), \ldots, \phi^i_{K_i}(\textbf{x}_j)$ for the approximation each $\frac{\partial u}{\partial \xi_i}, i=1, 2, \ldots, r$. Once the data-driven basis are constructed, we then use the collected data to train two neural networks to approximate the parameter-to-solution map, one for $\vxi \rightarrow \textbf{c}(\vxi)$ which gives $u(\textbf{x}_j;\vxi)=\sum_{k=1}^Kc_k(\vxi)\phi_k(\textbf{x}_j)$, and another one for $\vxi \rightarrow (\textbf{c}^1(\vxi), \textbf{c}^2(\vxi)\dots, \textbf{c}^r(\vxi))$ which gives $\frac{\partial u(\textbf{x}_j;\vxi)}{\partial \xi_i}=\sum_{k=1}^{K^i}c^i_k(\vxi)\phi^i_k(\textbf{x}_j), i=1, 2, \dots, r$, as described in Section~\ref{sec:AHMC}. In our experiments, the first network has $4$ hidden layers and $20$ units within each hidden layer. 
The second network has the same structure except there are $40$ hidden units in each hidden layer.

 We specify the number of leap-frog steps in \eqref{eq:Leapfrog} to be $10$, the step size $\Delta t=0.16$ for all methods. Typically, we start sampling from the posterior after observing mixing. For the standard HMC method and random network surrogate method \cite{zhang2017hamiltonian},  with which we compare, they share the same burn-in stage and starting point. We compute the relative error of the posterior mean up to a computation time $t$ by
\begin{equation}
\frac{\big|\big|\frac{1}{\#\{i:t_i\leq t\}}\sum_{i:t_i\leq t}\vxi_i - E_{\vxi|\textbf{y}}(\vxi)\big|\big|_2}{\big|\big|E_{\vxi|\textbf{y}}(\vxi)\big|\big|_2}.
\end{equation}

\begin{figure}
	\centering
	\begin{tabular}{cc}
		\includegraphics[width=0.48\linewidth]{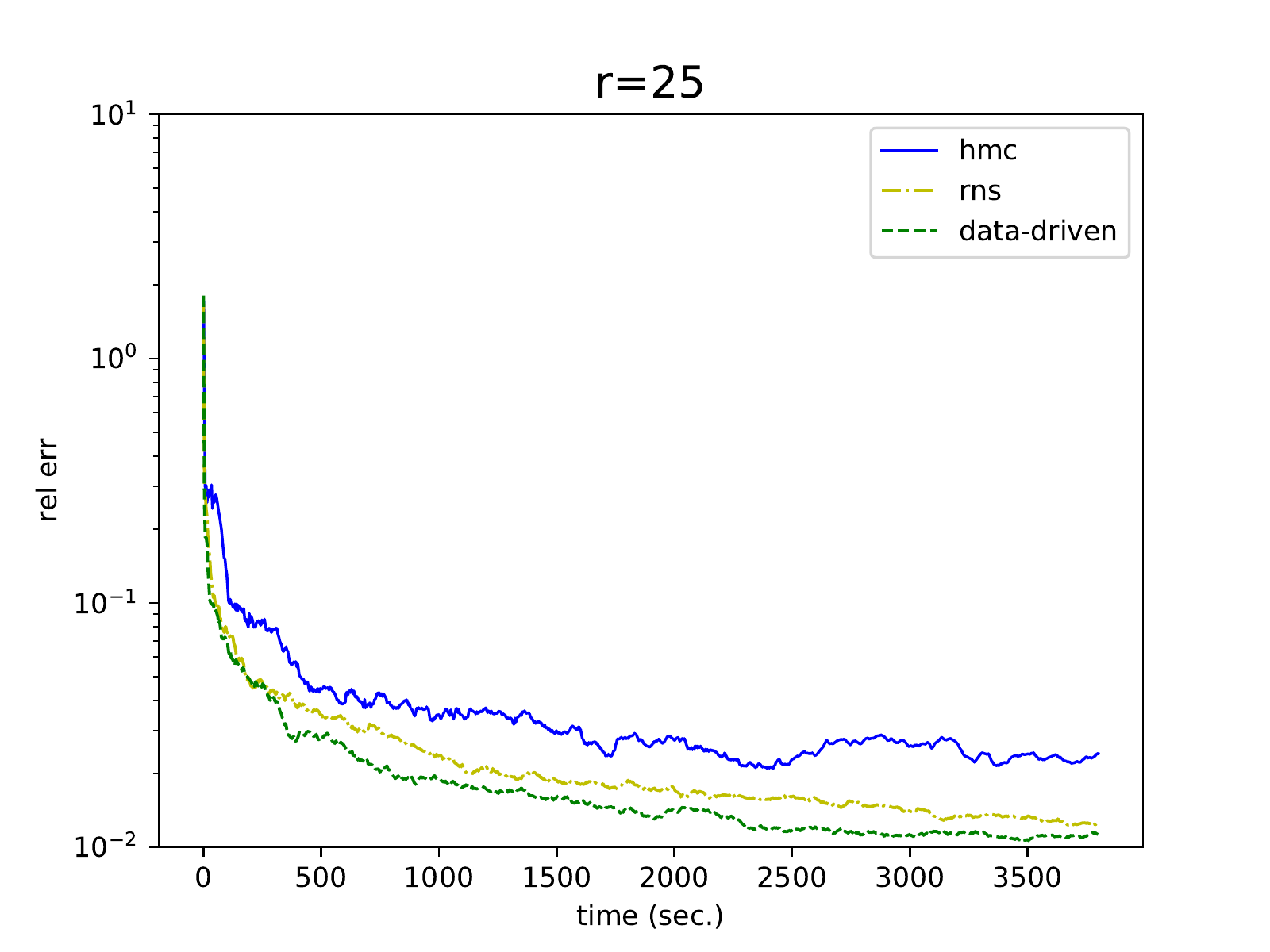} & \includegraphics[width=0.48\linewidth]{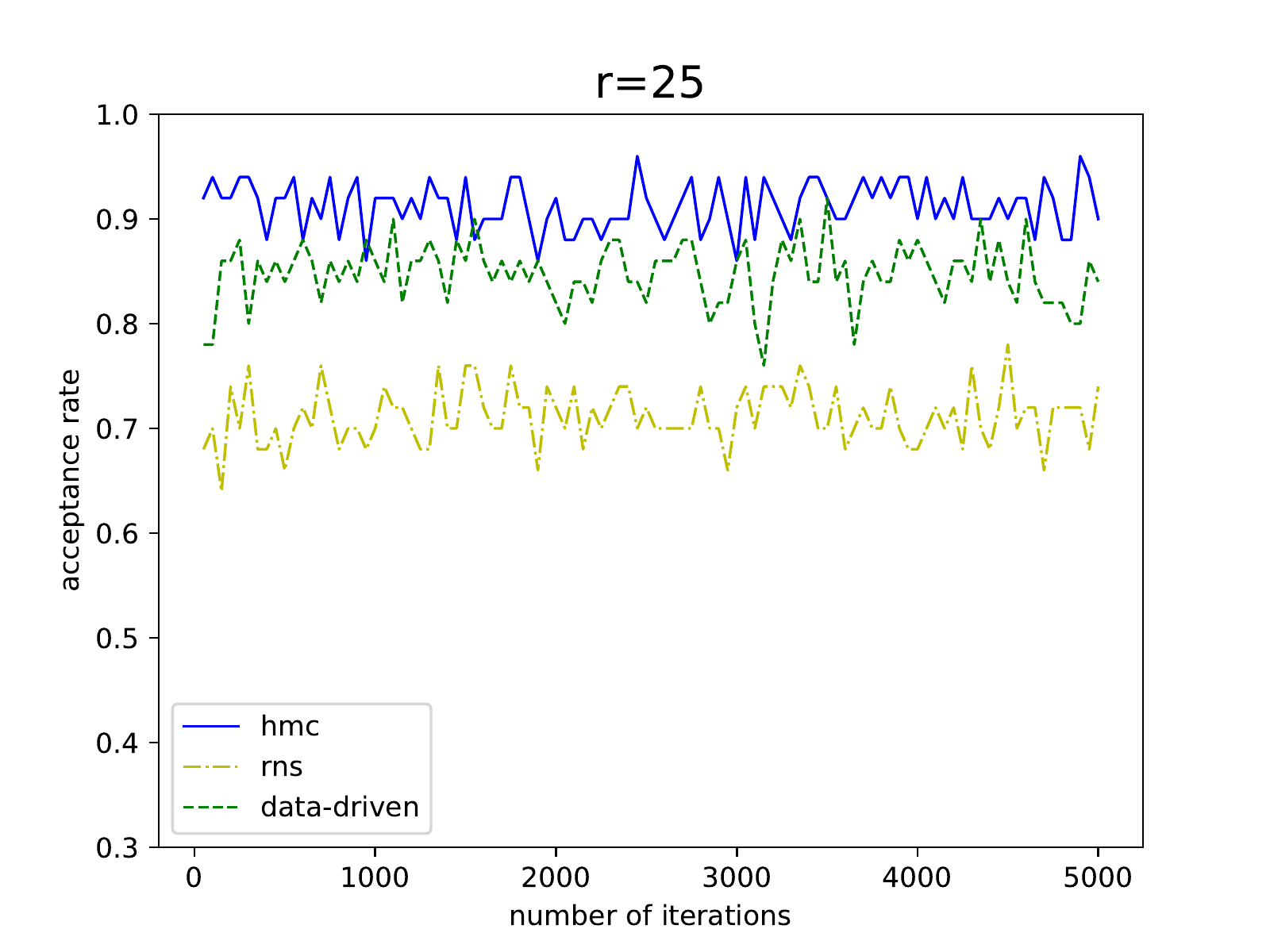}  \\
		\includegraphics[width=0.48\linewidth]{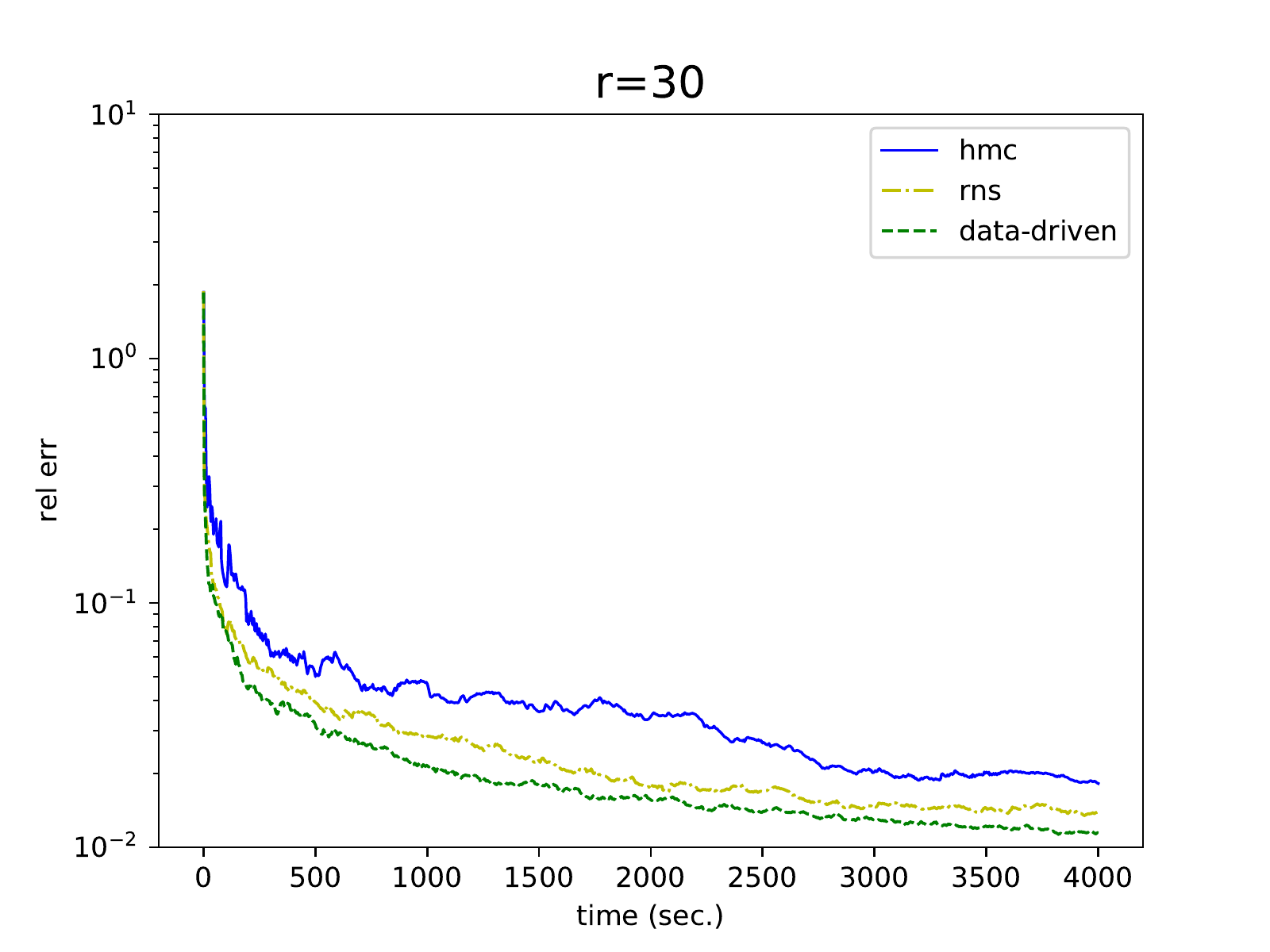} & \includegraphics[width=0.48\linewidth]{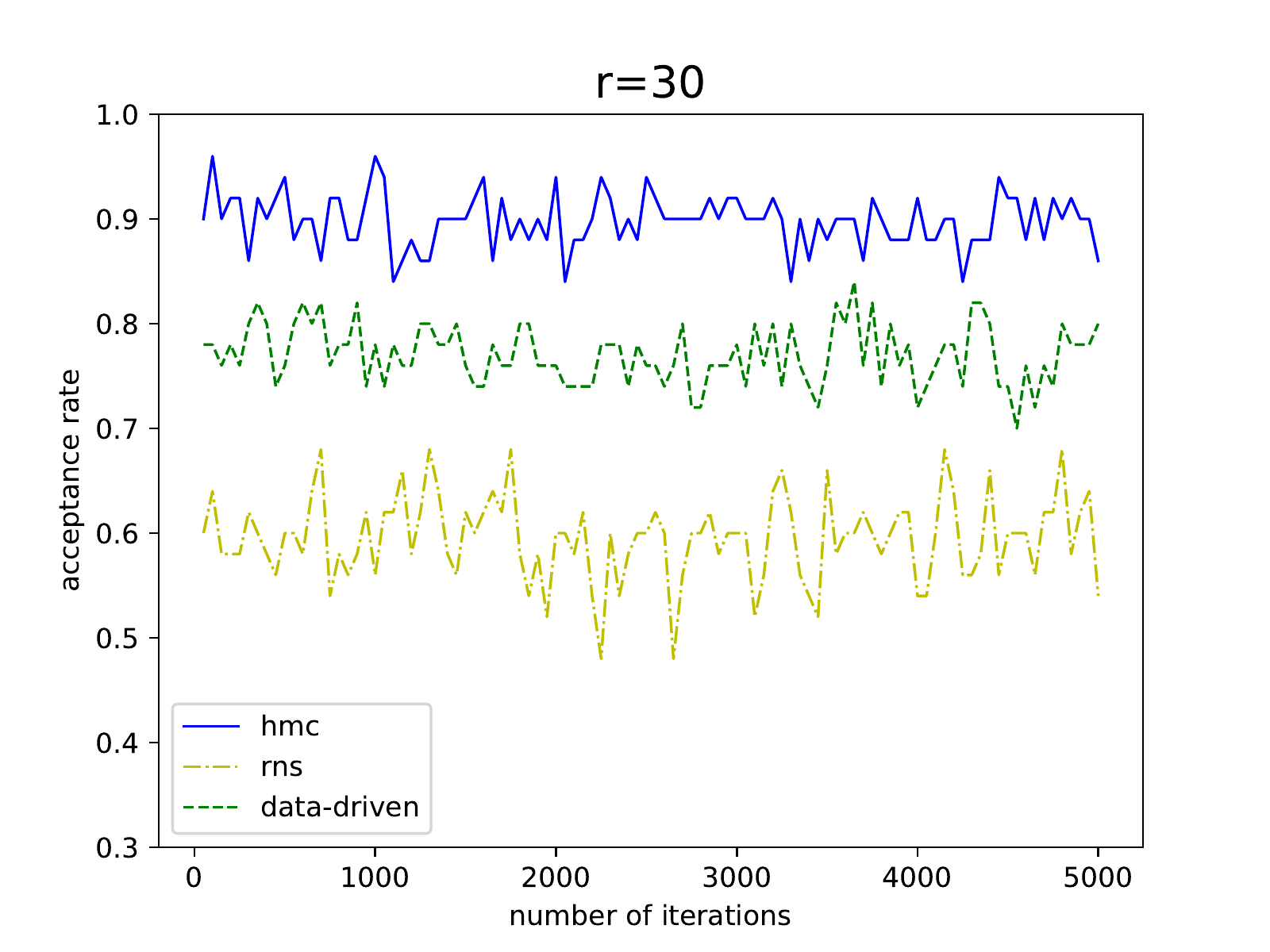} \\
		\includegraphics[width=0.48\linewidth]{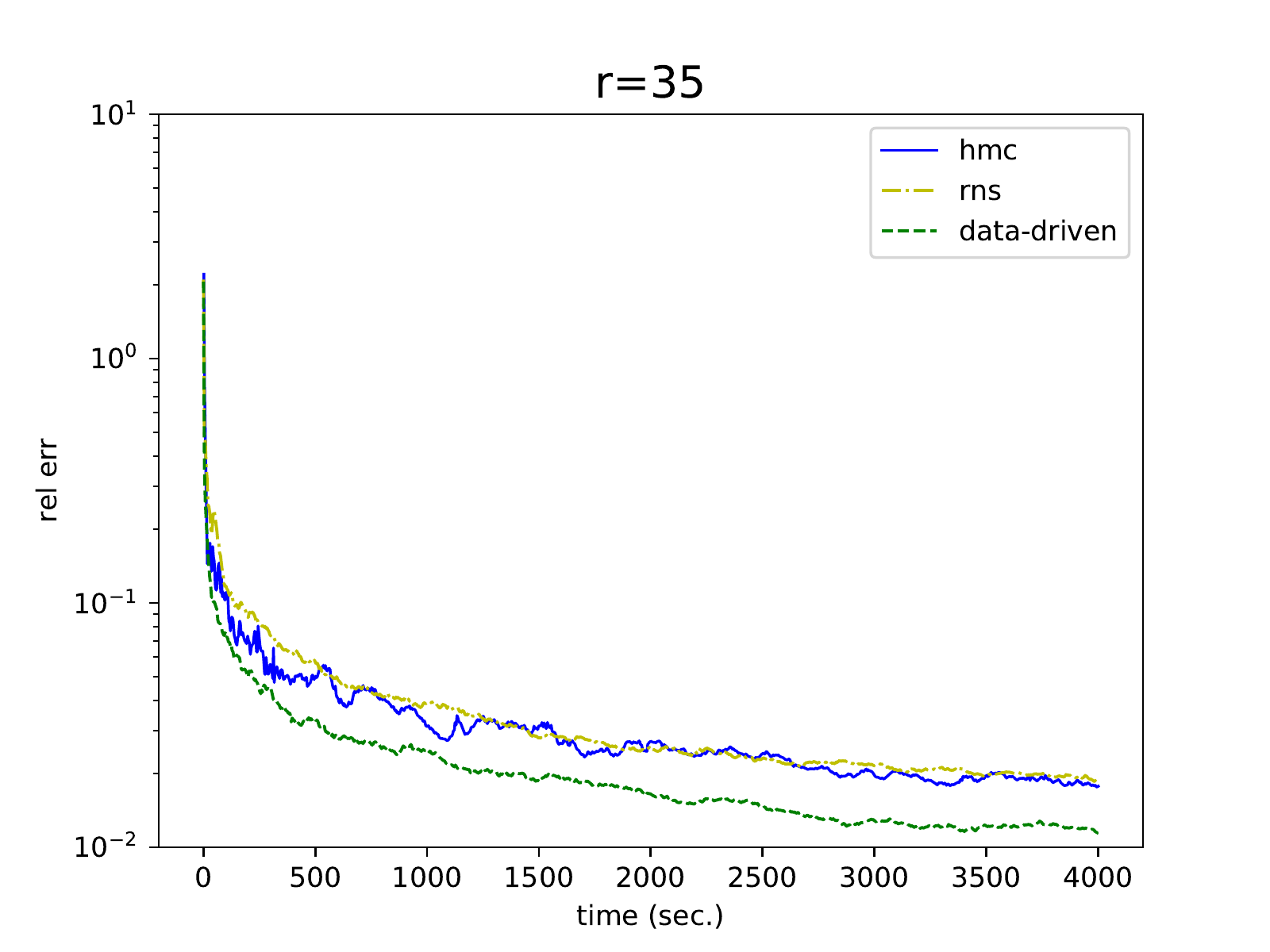} & \includegraphics[width=0.48\linewidth]{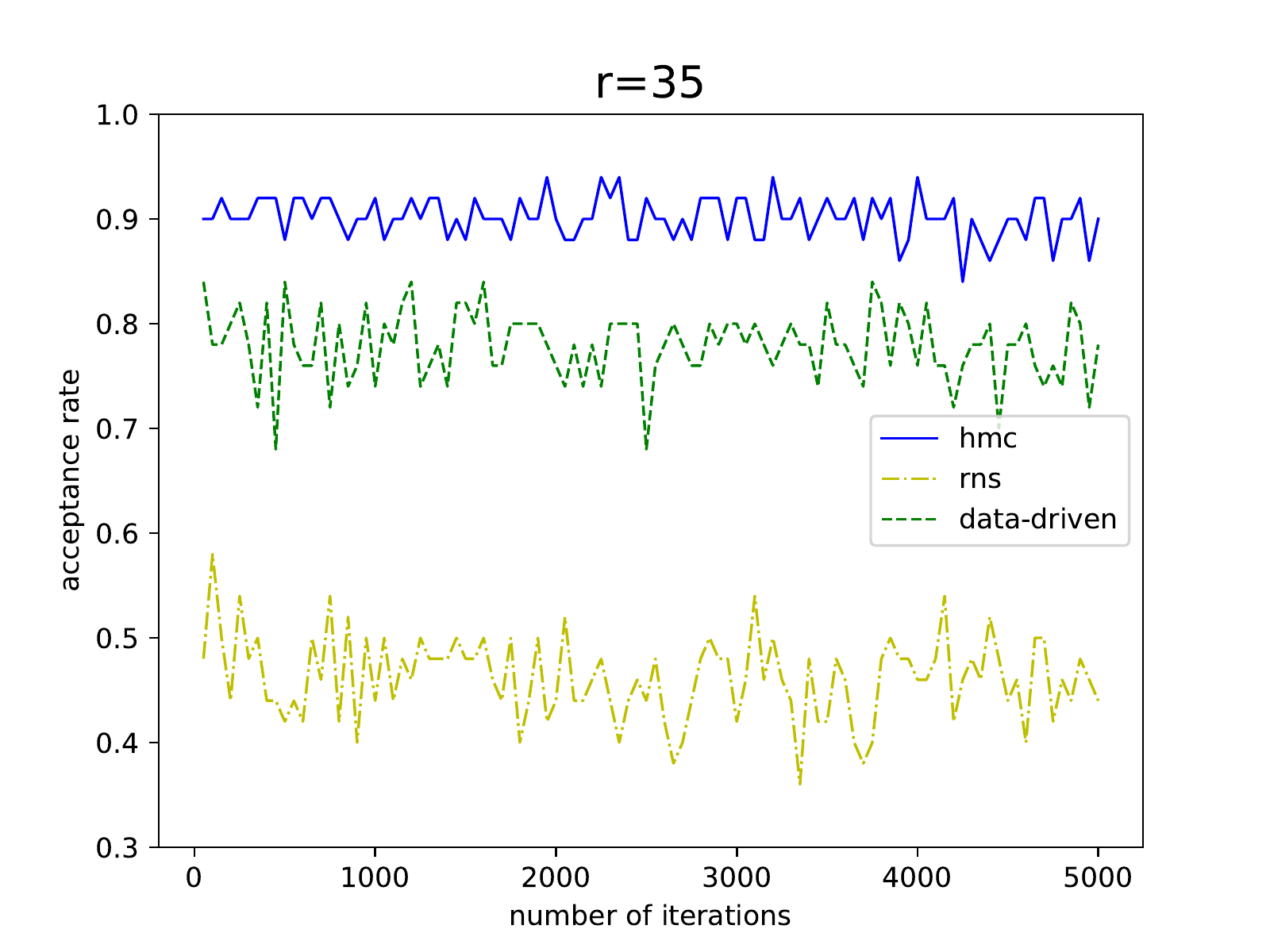} \\
		\multicolumn{1}{c}{relative errors w.r.t. computational time} & \multicolumn{1}{c}{acceptance rates w.r.t. iteration number} \\
		\hline
	\end{tabular}
	\caption{Numerical results for the random coefficient with $25$-, $30$- and $35$-dimensional inputs, where ``hmc'', ``rns'' and ``data-driven'' refer to the standard HMC method, the random network surrogate method and the accelerated HMC method, respectively.}\label{fig:error_and_acceptance}
\end{figure}

The left column of Figure \ref{fig:error_and_acceptance} plots the relative error of the posterior mean vs computation time in log scale for standard HMC, random network surrogate method, and the proposed accelerated HMC method for $r=25, 30, 35$. We see significantly improved performance of the proposed method.  The right column of Figure \ref{fig:error_and_acceptance} plots the corresponding acceptance rate for the proposal by Hamiltonian dynamics for each method. As we can see, the significant dimension reduction and the efficient neural network approximation of parameter-to-solution map does not compromise the acceptance rate much. Moreover, acceptance rate maintains high and stable as the input dimension increases for our model-based and data-driven approach.  For the random network surrogate method, the surrogate of the Hamiltonian in parameter space is based on least square approximation of sampled data, e.g, from burn-in stage, using a set of random basis. Since this approach is  purely data-driven without model knowledge, to maintain the approximation accuracy, the number of random basis has to increase with the dimension of the parameter space although the intrinsic dimension of the underlying model remains the same. In this experiment, we fix the number of random basis at 1000. As the input dimension increases, the approximation power of the surrogate using fixed number of basis decreases and hence the approximation error becomes larger and the acceptance rate drops quite sharply. 

Table \ref{tab:efficiency} shows the averaged time per each HMC iteration, averaged acceptance rate, effective sample size (min, median, max), and time normalized effective sample size for each method. The effective sample size is defined as 
\begin{equation}
ESS=B[1+2\sum_{k=1}^{K}\gamma(k)]^{-1}, 
\end{equation}
where $B$ is the number of MCMC samples and $\sum_{k=1}^{K}\gamma(k)$ is the sum of $K$ monotone sample autocorrelations. It shows that our model-based and data-driven approach has a good balance between the computation efficiency and exploration efficiency and hence achieves the best overall performance. 

\begin{table}[hbtp]
	\centering
	\begin{tabular}{|cc|ccccc|}
		\hline
		Dimension & Method & AR &  s/Iter & ESS & min(ESS)/s & med(ESS)/s \\
		\hline
		$r=25$ & hmc & 0.91 & 1.27 & (962 , 5000 , 5000) & 0.15 & 0.79 \\
		&rns  & 0.71 & 0.076 & (1866 , 2518 , 3095)  & 4.91 & 6.63 \\
		&data-driven & 0.85 & 0.094 & (3395 , 4307 , 4933)  &  7.22& 9.16 \\
		\hline
		$r=30$ & hmc & $0.90$ & 1.49 &  (1507 , 5000 , 5000) & 0.20   & 0.67 \\
		&rns  & $0.60$ & 0.082 &  (1286 , 1659 , 2069) &  3.14 & 4.05 \\
		&data-driven & $0.78$ & 0.108 &  (2905 , 3414 , 4116) &  5.38  &  6.32  \\
		\hline
		$r=35$ & hmc & 0.90 & 1.72 & (2061 , 5000 , 5000) & 0.24 & 0.58 \\
		&rns  & 0.47 & 0.095 & (597 , 968 , 1305) &  1.26 & 2.04\\
		&data-driven & 0.78 & 0.125&  (2823 , 3445 , 4030) & 4.52 & 5.51\\
		\hline
	\end{tabular}
	\caption{Comparisons of algorithms. The acceptance rate (AR), computational time for each iteration, effective sample size (ESS) and time-normalized ESS are provided.}\label{tab:efficiency}
\end{table}

\subsection{A log-normal coefficient with anisotropic heterogeneity}\label{sec:numericalex2}
\noindent
In the second example, a Gaussian prior with zero mean and covariance function
\begin{align}
c(\textbf{x},\textbf{x}^{\prime}) = \sigma_{a}^2\exp\Big(-\frac{|x_1-x_1^{\prime}|^2}{2l_1^2}-\frac{|x_2-x_2^{\prime}|^2}{2l_2^2}\Big)
\label{CovarianceFunction2}
\end{align}
is assumed on $\log(a(\textbf{x},\omega))$, where $\textbf{x}=(x_1, x_2)$ and $\textbf{x}^{\prime}=(x_1^{\prime}, x_2^{\prime})$ are any two points on $[0,1]\times[0,1]$, and $l_1$ and $l_2$ are the 
correlation lengths in $x_1$ and $x_2$. The diffusion coefficient is approximated via a truncated Karhunen-Lo\`eve (KL) expansion as in \eqref{KLE1}, only with a different prior covariance function \eqref{CovarianceFunction2}.

To generate the training data, we solve the elliptic problem \eqref{forwardex1} with the same boundary condition \eqref{forwardBC}. The discretization is done on a uniform grid with $65\times 65$ points through triangle finite element basis functions. We choose the number of truncated KL modes $r=30$, $\sigma=0.1$ be the noise in the observation data \eqref{observation1}, and $\sigma_{a}=0.5$, $l_1=0.08$ and $l_2=0.4$ in the prior covariance function \eqref{CovarianceFunction2}. All other settings are the same as in Section \ref{sec:numericalex1}. Suppose the measurements are placed on $17\times 17$ grids of the numerical solution $u(\textbf{x},\cdot)$, i.e. $m=289$ in \eqref{observation1}.  
We now infer the log-normal coefficient $\log(a(\textbf{x},\vxi))$ based on the observation data. 

\begin{figure}[hbtp]
	\centering
	\begin{subfigure}{0.49\textwidth}
		\includegraphics[width=1.1\linewidth]{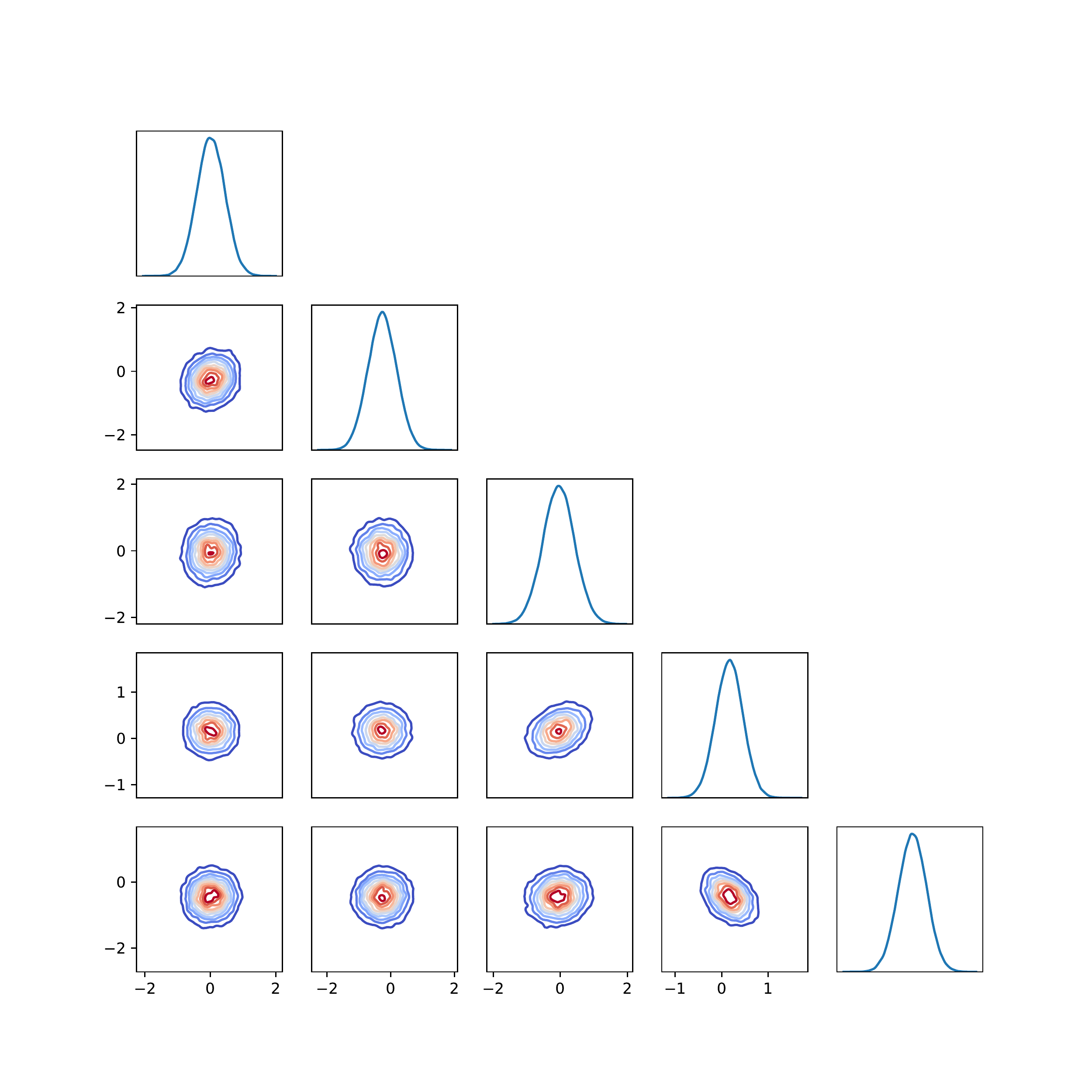}
		\subcaption{The standard HMC method.}
	\end{subfigure}
\centering
	\begin{subfigure}{0.49\textwidth}
		\includegraphics[width=1.1\linewidth]{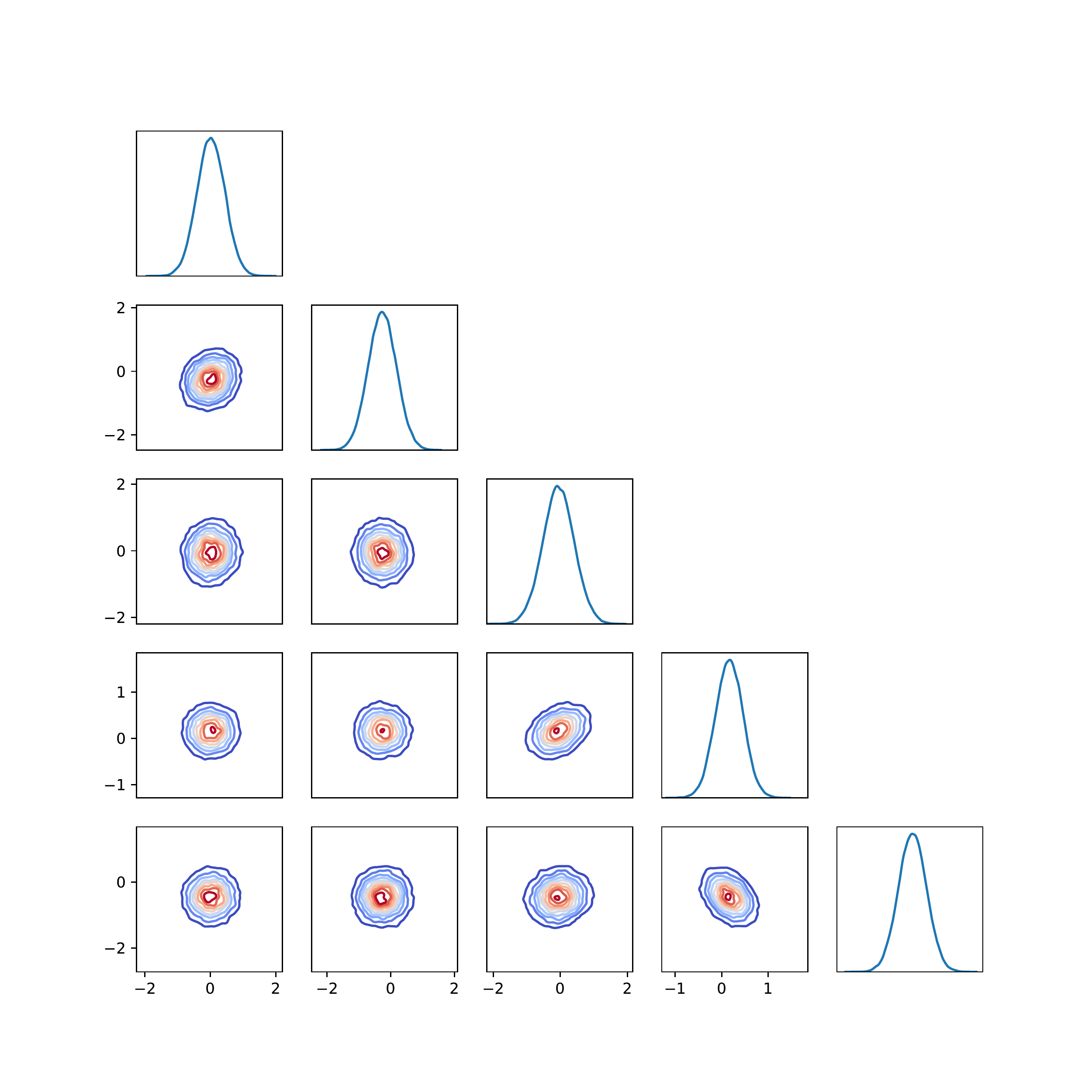}
		\subcaption{The accelerated HMC method.}
	\end{subfigure}
	\caption{Comparing one- and two-dimensional posterior marginals of $\xi_2$, $\xi_4$, $\xi_7$, $\xi_9$, $\xi_{13}$.}
	\label{fig:marginals}
\end{figure}
\begin{figure}[hbtp]
\centering
	\begin{tabular}{lc}

		\includegraphics[width=0.48\linewidth]{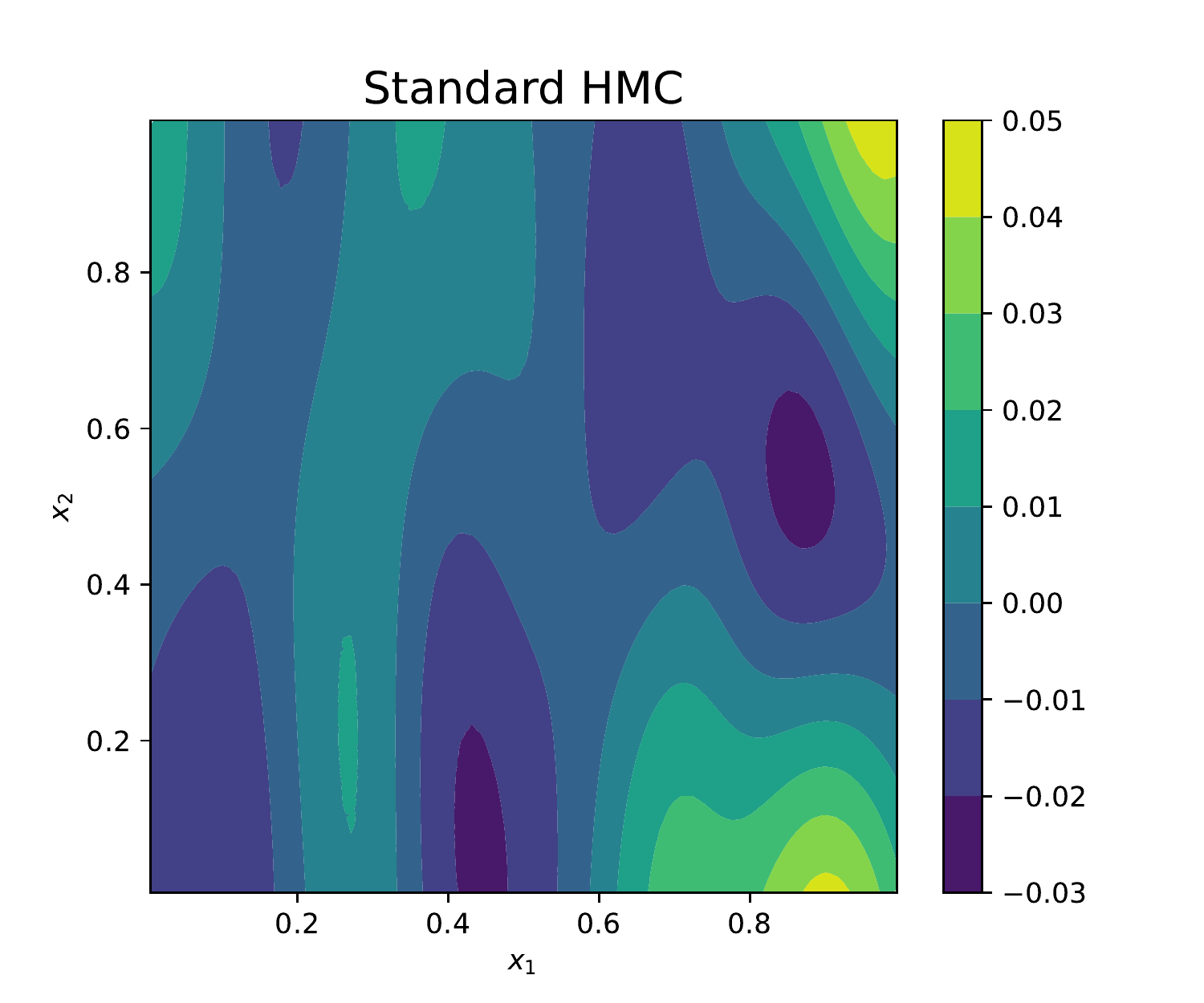} & 
		\includegraphics[width=0.48\linewidth]{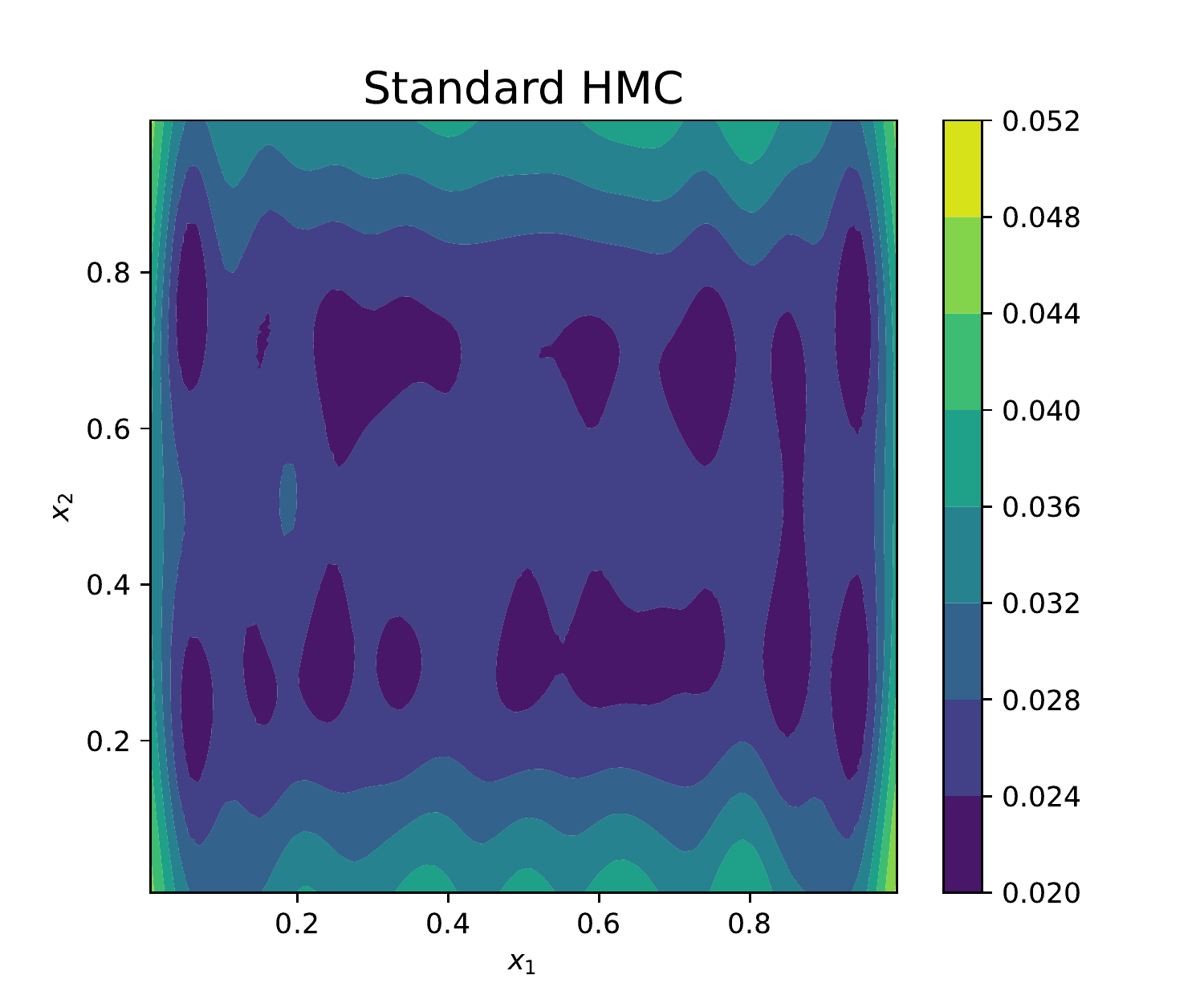} \\

		\includegraphics[width=0.48\linewidth]{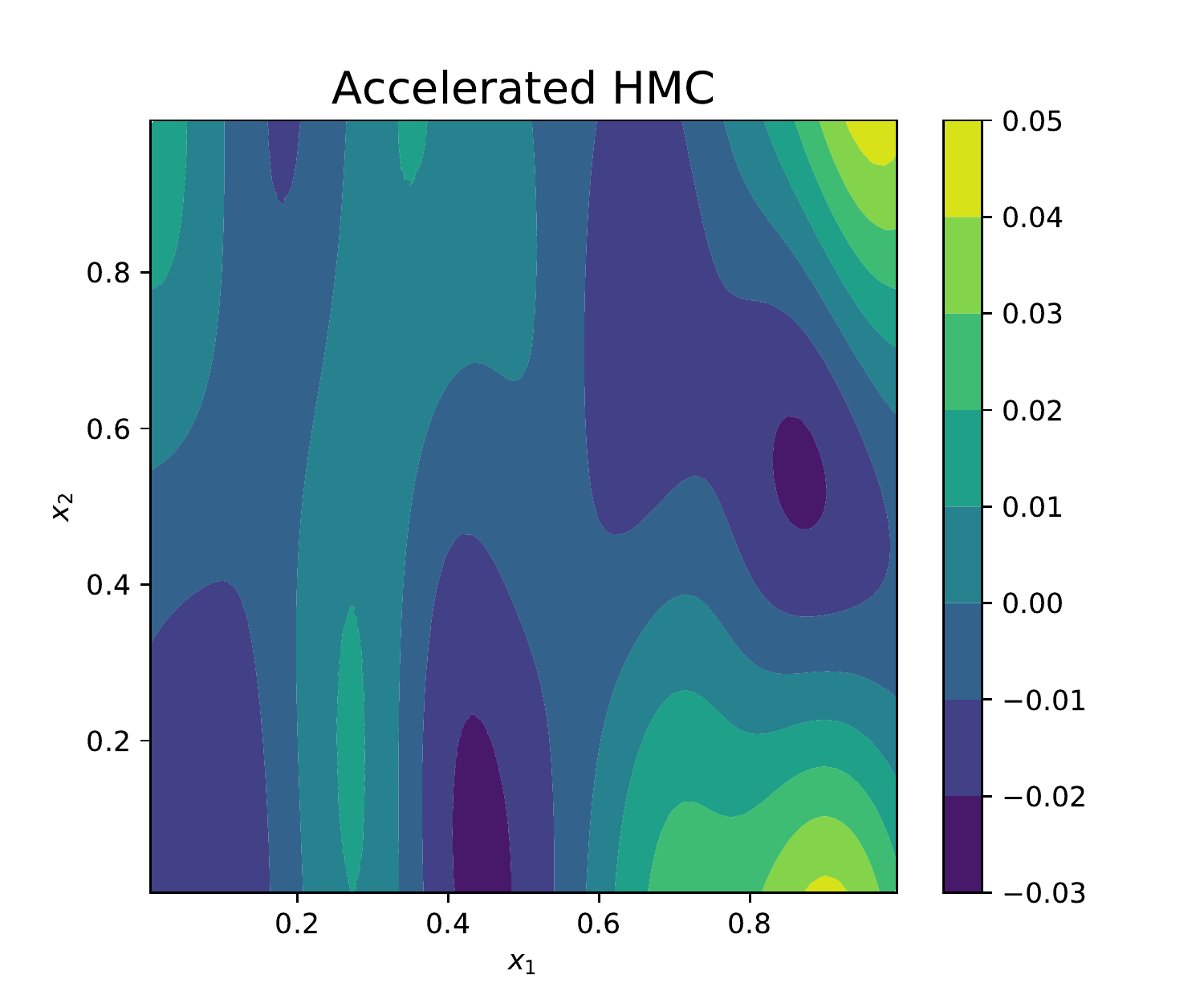} &
		\includegraphics[width=0.48\linewidth]{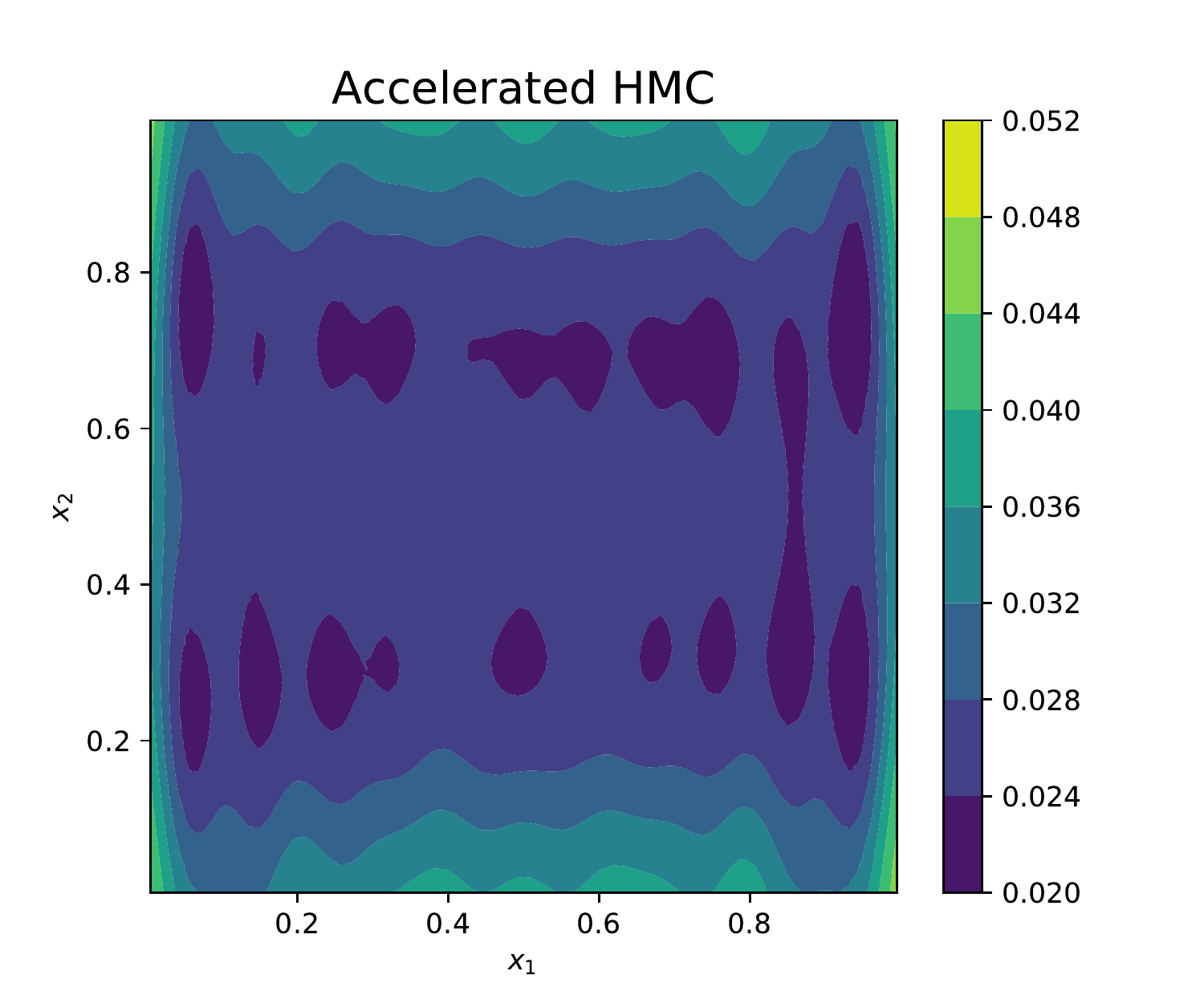} \\
		\multicolumn{1}{c}{Posterior mean.} & \multicolumn{1}{c}{Posterior standard deviation.}\\
	\end{tabular}
	\caption{Posterior statistics obtained by standard HMC and accelerated HMC.}\label{fig:ex2posthmc}
\end{figure}

To illustrate that our method indeed converges to the right target distribution, Figure \ref{fig:marginals} provides the one- and two- dimensional posterior marginals of some selected parameters obtained by standard HMC and the accelerated HMC.
Figure \ref{fig:ex2posthmc} shows the posterior mean and posterior standard deviation obtained by the standard HMC method and  the accelerated HMC method, respectively. The relative errors of the posterior mean and posterior standard deviation are $0.047$ and $0.024$.
Therefore, with the accelerated HMC method, we can significantly reduce the computation cost (by almost an order of magnitude in this case) while maintaining the approximation accuracy of the standard HMC.

\begin{figure}[hbtp]
\centering
\begin{subfigure}{0.48\textwidth}
\includegraphics[width=1.0\linewidth]{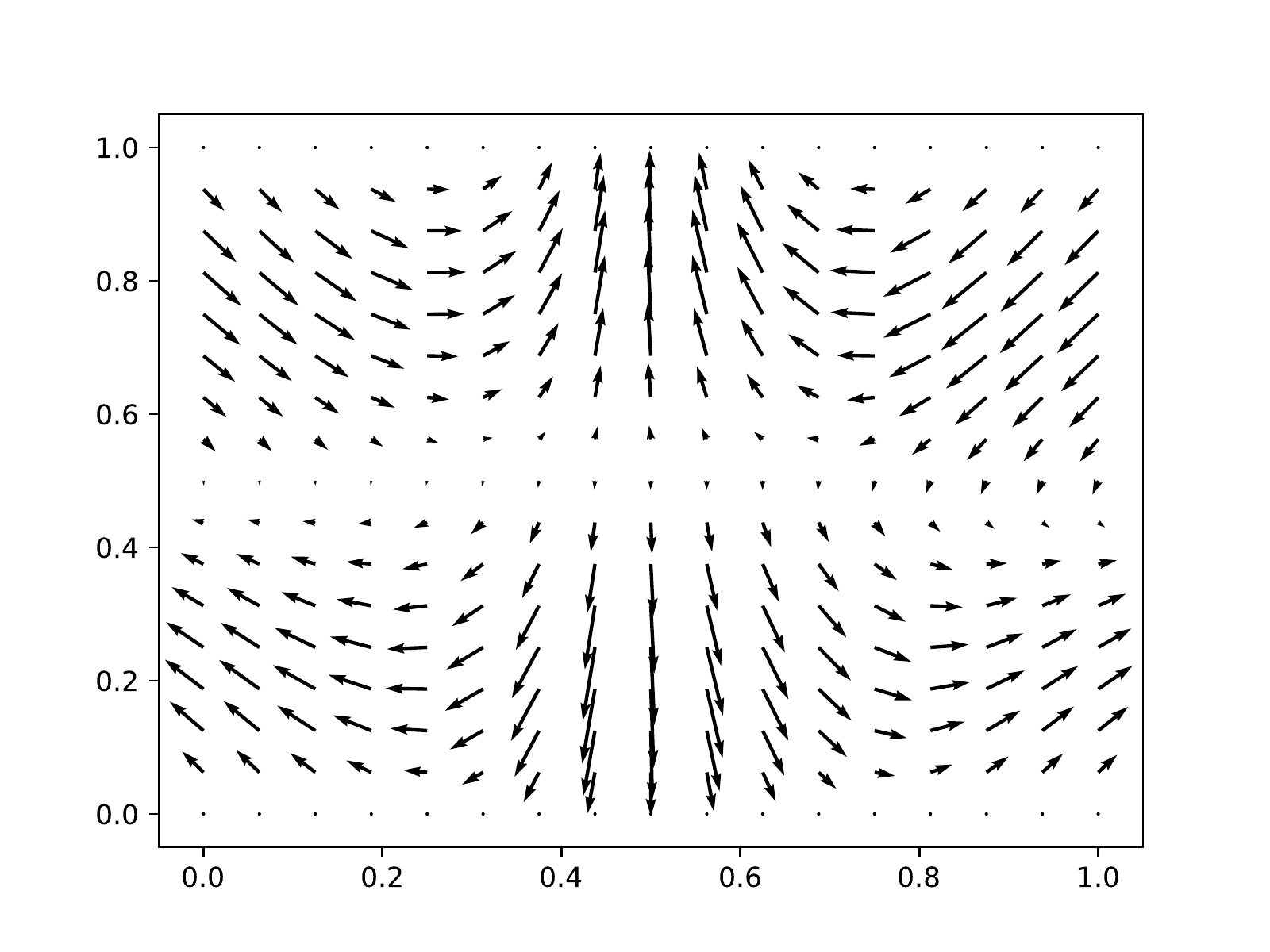}
\subcaption{The standard HMC method.}
\label{fig:ex2partiald-hmc}
\end{subfigure}
\begin{subfigure}{0.48\textwidth}
\includegraphics[width=1.0\linewidth]{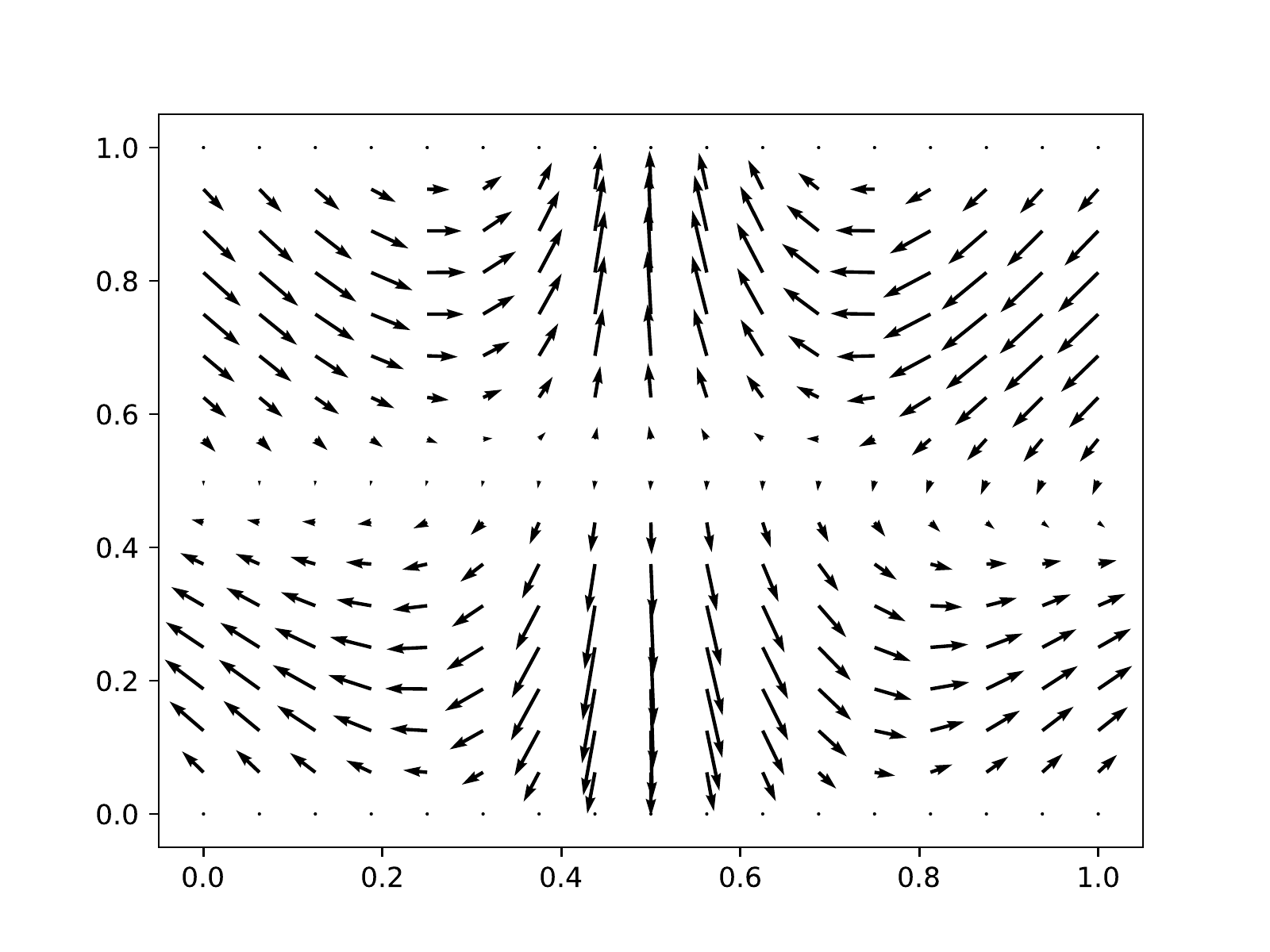}
\subcaption{The accelerated HMC method.}
\label{fig:ex2partiald-dd}
\end{subfigure}
\caption{Partial derivative of $u(\textbf{x},\vxi)$ with respect to $\xi_1$ and $\xi_2$.}
\label{fig:PartialDerivative}
\end{figure}

Finally, we compare the partial derivatives $\frac{\partial u(\textbf{x},\vxi)}{\partial \xi_1}$ and $\frac{\partial u(\textbf{x},\vxi)}{\partial \xi_2}$ at the approximate MAP state obtained via standard HMC method and the accelerated HMC method in Figure \ref{fig:PartialDerivative}. The relative errors of $\frac{\partial u(\textbf{x},\vxi)}{\partial \xi_1}$ and $\frac{\partial u(\textbf{x},\vxi)}{\partial \xi_2}$ are $0.013$ and $0.019$, respectively. We also examine the relative errors at ten posterior sample of $\vxi$ and the result is presented in Table \ref{table:temp1}. These results demonstrate the effectiveness of our intrinsic low dimensional data-driven basis on providing fast and accurate gradient approximations for accelerating Hamiltonian dynamics.



\begin{table}[h]
	\centering
	\begin{tabular}{|ccccccccccc|}
		\hline
		$\frac{\partial u(\textbf{x},\vxi)}{\partial \xi_1}$ &  0.032 & 0.014& 0.036& 0.054& 0.016& 0.021& 0.032 & 0.103& 0.028 & 0.018 \\
		\hline
		$\frac{\partial u(\textbf{x},\vxi)}{\partial \xi_2}$ & 0.064  & 0.055& 0.022& 0.045& 0.024& 0.037& 0.043 & 0.067& 0.080 & 0.032
		\\      
		\hline
	\end{tabular}
	\caption{relative errors between the exact solution and the data-driven approximation at ten random samples from the posterior of $\vxi$.}
	\label{table:temp1}
\end{table}

\section{Conclusion} \label{sec:conclusion}
\noindent 
The HMC method can generate less correlated proposals with high acceptance probabilities, which greatly improves the performance of the MCMC methods in solving Bayesian inverse problems. However, when applying the HMC method to solve a Bayesian inverse problem modeled by elliptic partial differential equations, one needs to compute solution to the elliptic PDEs and their derivatives repeatedly in order to generate data and evaluate the Hamiltonian, which makes the HMC method extremely expensive. 

By exploiting the intrinsic low-dimensional structures of the underlying model and constructing a data-driven basis, our proposed method achieves significant dimension reduction in the solution space. Then, equipped with the data-driven basis, neural networks are trained as efficient approximations of the parameter-to-solution maps, which significantly reduce the computation cost in obtaining the PDE solution and its derivatives for the Hamiltonian dynamics in proposing a new sample. Through numerical tests, we demonstrate that our method strikes a good balance between computation efficiency and exploration efficiency and provides an effective data and model-based approach for elliptic Bayesian inverse problems.



\section*{Acknowledgements}
\noindent
The research of S. Li is partially supported by the Doris Chen Postgraduate Scholarship. 
The research of C. Zhang is supported by the Key Laboratory of Mathematics and Its Applications (LMAM)
and the Key Laboratory of Mathematical Economics and Quantitative Finance (LMEQF) of Peking
University.  The research of Z. Zhang is supported by the Hong Kong RGC General Research Fund Project 17300318, Seed Funding Programme for Basic Research (HKU), and Basic Research Programme (JCYJ20180307151603959) of The Science, Technology, and Innovation Commission of Shenzhen Municipality. The research of H. Zhao is partially supported by NSF grant DMS-2048877 and DMS-2012860. 
The computations were performed using research computing facilities offered by Information Technology Services, the University of Hong Kong. 

\bibliographystyle{plain}
\bibliography{ZWpaper}


\end{document}